\documentclass[11pt,a4paper]{amsart}
\usepackage[top=30mm,bottom=30mm,left=27mm,right=27mm]{geometry}
\usepackage{amssymb,amsthm,amsmath,mathtools}
\usepackage{enumerate,verbatim,microtype,needspace}
\usepackage[colorlinks=true,citecolor=blue]{hyperref}
\allowdisplaybreaks[1]
\DeclareMathAlphabet{\mathcal}{OMS}{cmsy}{m}{n}
\SetMathAlphabet{\mathcal}{bold}{OMS}{cmsy}{b}{n}

\newtheorem{thmx}{Theorem}

\newtheorem{theorem}{Theorem}[section]
\newtheorem{proposition}[theorem]{Proposition}
\newtheorem{lemma}[theorem]{Lemma}
\newtheorem{corollary}[theorem]{Corollary}
\theoremstyle{definition}
\newtheorem{definition}[theorem]{Definition}

\newtheorem{remark}[theorem]{Remark}
\numberwithin{equation}{section}

\newcommand{\N}{\mathbb N}
\newcommand{\Z}{\mathbb Z}
\newcommand{\Q}{\mathbb Q}
\newcommand{\R}{\mathbb R}
\newcommand{\C}{\mathbb C}
\newcommand{\T}{\mathbb T}
\newcommand{\E}{\mathbb E}
\newcommand{\I}{\mathcal I}
\newcommand{\HK}{\mathcal Z}
\newcommand{\W}{\mathcal W}
\newcommand{\Krat}{\mathcal K_{\mathrm{rat}}}
\newcommand{\roots}{\mu_{\infty}}
\newcommand{\av}[2]{\frac1{#1}\sum_{#2=1}^{#1}}
\newcommand{\norm}[1]{\left\lVert#1\right\rVert}

\DeclareMathOperator{\supp}{supp}
\DeclareMathOperator{\id}{id}

\begin{document}
	\title[Phase factors and optimal seminorms for polynomial ergodic averages]{Phase factors and optimal seminorms for polynomial averages of commuting transformations}
	\author{Zhuowei Liu}
	\date{\today}
	\address[Z. Liu]{School of Mathematics (Zhuhai), Sun Yat-sen University,
		Zhuhai, Guangdong, 519000, P.R. China}
	\email{liuzhw55@mail2.sysu.edu.cn}
	
\begin{abstract}
For invertible commuting measure-preserving transformations $T_1, \ldots, T_\ell$ on a standard probability space, we study polynomial averages along both integers and primes. We provide a rank-based criterion ensuring that, even when the polynomials are linearly dependent, finite-order phase factors are characteristic factors. For polynomials with integer coefficients and zero constant terms, it suffices that the dimension of the space spanned (over $\mathbb{Q}$) by any nonempty subfamily of $k$ polynomials exceeds $k/2$. These phase factors are generated by polynomial phases associated with positive integer powers of the individual transformations and may be strictly smaller than the corresponding Host-Kra factors of the same order. The required order depends only on the number of polynomials and their maximum degree. A four-term example involving a totally ergodic Heisenberg nilsystem shows that this criterion may fail if the equality $k/2$ is permitted.

For a family of three polynomials, this criterion is equivalent to the condition that the polynomials are pairwise linearly independent. Consequently, it applies to the family $t, t^2, t^2+t$, where the third polynomial is the sum of the first two. For this specific family, we prove that the optimal orders for the phase factors and the Host-Kra factors are $4, 2, 2$, respectively--a result that holds for both integers and primes. The corresponding optimal seminorms are $U^5(T_1)$, $U^3(T_2)$, and $U^3(T_3)$, which provides an answer to a question posed by Kuca.
\end{abstract}
	\maketitle

	\section{Introduction}
	\subsection{Background and main questions}
	The study of polynomial ergodic averages originated from the connection between multiple recurrence and arithmetic configurations. Furstenberg’s ergodic proof of Szemer\'edi's theorem \cite{Furstenberg1977} connected multiple recurrence with arithmetic progressions. The polynomial recurrence theorem of Bergelson and Leibman \cite{BergelsonLeibman1996} extended this connection to
	polynomial configurations and commuting transformations. 
	A central problem is to characterize the parts of a dynamical system
	that determine the limits of these averages.
	
	Let $\ell\ge3$, $[\ell]:=\{1,\ldots,\ell\}$, and let
	$T_1,\ldots,T_\ell$ be invertible commuting measure-preserving
	transformations of a standard probability space $(X,\mathcal X,\mu)$.
	For integer polynomials $\mathbf p=(p_1,\ldots,p_\ell)$ with
	$p_i(0)=0$, write
	\begin{equation}\label{eq:general-average}
		A_N^{\mathbf p}(f_1,\ldots,f_\ell)
		=\frac1N\sum_{n=1}^N\prod_{i=1}^\ell T_i^{p_i(n)}f_i,
		\qquad f_i\in L^\infty(\mu),
	\end{equation}
	where $Tf=f\circ T$.
	For linear and polynomial iterations, Tao \cite{Tao2008} and Walsh \cite{Walsh2012} proved norm convergence, respectively.
	To describe limits more precisely, one seeks a collection of factors $(\mathcal A_i)_{i=1}^\ell$ such that, upon replacing each $f_i$ with 	$\E(f_i\mid\mathcal A_i)$, the difference in the average tends to zero in $L^2(\mu)$. Such a collection of factors is called a system of characteristic factors.
	
	For a single transformation, the Host–Kra–Ziegler structure theory \cite{HostKra2005,Ziegler2007}, combined with the theory of equidistribution of polynomials on nilmanifolds, provides the characteristic factors for polynomial averages\cite{HostKraPolynomial2005,LeibmanConvergence2005}.
	
	For several commuting transformations, people wish to use the factors corresponding to each $T_i$ separately.
	Chu, Frantzikinakis, and Host established Host–Kra factor control for various transformations, applicable to polynomials with distinct degrees
	\cite{ChuFrantzikinakisHost2011}.
	Following the estimates on the seminorm in \cite{DonosoFerreMoraguesKoutsogiannisSun2024, DonosoKoutsogiannisSun2022}, Frantzikinakis and Kuca proved that for any family of pairwise linearly independent polynomials, the finite-order Host–Kra factors of the individual transformations form a system of characteristic factors \cite[Theorem~2.8]{FrantzikinakisKuca2025}. This theorem does not require ergodicity, and the necessary order depends only on the number of polynomials and their maximum degree. Section 2.4 of the same paper also discusses a version of this result for prime numbers. For a family of linearly independent polynomials, Frantzikinakis and Kra proved that the rational Kronecker factors are characteristic for averages taken over integers for a single ergodic transformation \cite[Theorem~1.1]{FrantzikinakisKra2006}. Frantzikinakis and Kuca extended this result to commuting transformations, covering both averages over integers and averages over prime numbers \cite[Theorems~2.10 and~2.13]{FrantzikinakisKuca2025}.
	
	For a family of polynomials that are pairwise linearly independent but collectively satisfy a linear relation, the above results leave two questions open: what is the required degree for each coordinate, and which parts of the corresponding Host–Kra factor determine the limit of the average? The first of these questions already arises in the following average:
	\begin{equation}\label{eq:average}
		A_N(f_1,f_2,f_3)
		=\frac1N\sum_{n=1}^N T_1^n f_1\,T_2^{n^2}f_2\,T_3^{n^2+n}f_3.
	\end{equation}
	The family of polynomials $t,t^2,t^2+t$  are pairwise linearly independent, but the third polynomial is the sum of the first two.
	When $T_1=T_2=T_3=T$ is ergodic, 
	the Kronecker factor $\HK_1(T)$ coincides with the characteristic factor \cite[Theorem~B]{Frantzikinakis2008}.
	For general commuting transformations, the existence of a control by a factor of finite order has been established \cite[Theorem~2.8]{FrantzikinakisKuca2025}.
	In \cite[Section~6.3, Problem~36]{Kuca2026}, Kuca asks about the optimal order, noting that it is not even clear whether $\HK_1(T_1),\HK_1(T_2),\HK_1(T_3)$  are sufficient.
	This problem examines whether the low-order description in a single transformation case can be extended to several commuting transformations.
	
	The second question is inspired by existing work on the single-transformation case, which reduces the problem to affine systems and factors generated by quasi-eigenfunctions \cite{Frantzikinakis2008}.
	The study of Weyl systems and polynomial complexity
	\cite{BergelsonLeibmanLesigne2007,Kuca2023Complexity}
	relates characteristic orders to affine models and polynomial
	identities. Even if the Host-Kra complexity is equal to the Weyl complexity, it cannot be concluded that, in a given commutative system, each Host-Kra factor can be replaced with a subfactor generated by polynomial phases for positive powers of the corresponding transformation.
	We investigate when this substitution holds and allow for linear relationships within a family of polynomials.
	
	In this paper, we answer the first question for \eqref{eq:average}: for all commuting systems, the optimal coordinatewise Host--Kra
	factor orders are $4,2,2$. Therefore,
	The first-order conclusion for a single transformation no longer holds in general commuting systems.
	For a pair of commuting affine transformations  on $\T^5$, one obtains three matching lower bounds simultaneously, with the possibility of setting $T_2 = T_3$.

    The upper bound comes from a more general reduction result, namely, reduction to a factor generated by polynomial phases for positive powers of the corresponding transformation.
	We prove that as long as the dimension of the space spanned by each nonempty subfamily of polynomials over the rational number field is greater than half its size, these smaller factors constitute a set of characteristic factors. This criterion covers all triples of pairwise linearly independent polynomials and also applies to longer families of polynomials with linear relationships.
	This reduction and the  order $4,2,2$ hold true for both integer and prime averages.
	Furthermore, the orders of the phase factors for this three-term polynomial family are optimal.

	\subsection{Phase factors and a rank criterion}
	
	Now we define the phase factor and state the result precisely.
	Let $\mathbb P$ be the primes and $\pi(N)=|\mathbb P\cap[1,N]|$.
	For $N\ge2$, put
	\begin{equation}\label{eq:prime-average}
		A_{N,\mathbb P}^{\mathbf p}(f_1,\ldots,f_\ell)
		=\frac1{\pi(N)}\sum_{\substack{p\le N\\p\in\mathbb P}}
		\prod_{i=1}^\ell T_i^{p_i(p)}f_i.
	\end{equation}
	Let $A_{N,\mathbb P}$ denote the prime-number version of \eqref{eq:average}. For averages over primes, the definition of the characteristic factor tuple remains the same. By Walsh’s theorem and the comparison results in \cite{FrantzikinakisHostKra2013}, these averages over primes converge.
	
	We adopt the following conventions: $\|f\|_{U^{d+1}(T)}=0$ if and only if
	$\E(f\mid\HK_d(T))=0$; see \cite{HostKra2005,HostKra2018}.
	For a function $\phi$ of modulus one, let
	$\Delta_T\phi=(T\phi)\overline\phi$, and define its iterates recursively. For $d\ge1$, let
	\begin{equation*}\label{eq:Wdefinition}
		\W_d(T)=\sigma_\mu\left(
		\left\{\phi\in L^\infty(\mu):|\phi|=1,\quad
		\Delta_{T^q}^{d+1}\phi=1\text{ for some }q\ge1\right\}\right),
	\end{equation*}
	where $\sigma_\mu$ denotes the completed sigma-algebra generated by
	the indicated functions. We call $\W_d(T)$ the \emph{phase factor}
	of order $d$.
	
	The functions satisfying $\Delta_{T^q}^{d+1}\phi=1$ are classical quasi-eigenfunctions for positive integer powers of $T$. Abramov established a structural theory for totally ergodic systems with quasi-discrete spectrum \cite{Abramov1962};
	see also \cite{HaaseMoriakov2018}.
	For  a fixed $q$, following the convention in \cite[Section~2.3]{HaaseMoriakov2018}, these generating functions are of order at most $d$ with respect to $T^q$.
	Equivalently, adopting the notation from
	\cite[Definition~1.13]{BergelsonTaoZiegler2010},
	$\W_d(T)=\bigvee_{q\ge1}\operatorname{Abr}_{<d+1}(T^q)$,
	where the completion is taken. The equivalence of different definitions of quasi-eigenfunction will be reviewed in Section~\ref{sec:preliminaries}. It is important to consider positive integer powers here because, even if a transformation is ergodic, the space of higher-order quasi-eigenfunctions may change upon taking powers
	\cite[Theorem~1.4(i)]{Frantzikinakis2006WW}.
	
	Lemma~\ref{lem:phasefactors} establishes that $\W_d(T)\subseteq\HK_d(T)$,
	and this inclusion can be strict, even when considering all positive integer powers.
	For an ergodic translation on the standard Heisenberg nilmanifold,
	$\W_2(T)=\HK_1(T)\subsetneq\HK_2(T)$; see
	\cite[Propositions~2.1(i), 2.4 and Theorem~2.2]{Frantzikinakis2008}
	and Section~\ref{subsec:heisenberg-application}.
	Therefore, even if the corresponding Host–Kra seminorm is not zero, the phase factor theorem may still deduce that the average tends to zero.
	
	For $J\subseteq[\ell]$, put
	$r_{\mathbf p}(J)=\dim_{\Q}\operatorname{span}_{\Q}\{p_j:j\in J\}$.
	The following criterion can be verified using the coefficient matrices of the respective families of polynomials.
	
	\begin{thmx}
		\label{thm:phase-reduction}\label{prop:reduction}
		Let $D\ge1$ and $\ell\ge3$. There exists an integer $s=s(D,\ell)\ge1$
		with the following property. Let $p_1,\ldots,p_\ell\in\Z[t]$ have
		degree at most $D$ and satisfy $p_i(0)=0$. Assume that
		\begin{equation}\label{eq:rank-condition}
			|J|<2r_{\mathbf p}(J)
			\qquad(\varnothing\ne J\subseteq[\ell]).
		\end{equation}
		For every system of invertible commuting measure-preserving
		transformations $T_1,\ldots,T_\ell$, the factors
		$\W_s(T_1),\ldots,\W_s(T_\ell)$ are characteristic for both
		\eqref{eq:general-average} and \eqref{eq:prime-average}.
		More precisely, if $f_i\in L^\infty(\mu)$ and
		$\E(f_j\mid\W_s(T_j))=0$ for some $j\in[\ell]$, then
		\[
		\|A_N^{\mathbf p}(f_1,\ldots,f_\ell)\|_2\longrightarrow0,
		\qquad
		\|A_{N,\mathbb P}^{\mathbf p}(f_1,\ldots,f_\ell)\|_2\longrightarrow0.
		\]
	\end{thmx}
	
	The rank condition excludes the zero polynomial and implies pairwise linear independence. For three polynomials, it is exactly equivalent to each being pairwise linearly independent. The order $s(D,\ell)$ is independent of the
	coefficients, the system,  and the functions involved.
	
	At the boundary where equality holds, there is a specific obstruction to the strict inequality.
	Proposition~\ref{prop:rank-boundary} examines
	$p_j(t)=t+jt^2$, $0\le j\le3$.
	Every nonempty subfamily satisfies $|J|\le2r_{\mathbf p}(J)$,
	with equality occurring only for the full family of four polynomials.
	On a totally ergodic Heisenberg nilsystem, the same set of four functions are orthogonal to $\W_s(T)$ for every $s\ge1$, while their integer average has a non-zero limit. Here, all four transformations are identical. Therefore, increasing the phase order cannot remove this obstacle.
	This example excludes the possibility of uniformly replacing $<$ by $\le$ in
	\eqref{eq:rank-condition}.
	
	By Lemma~\ref{lem:rank-cover}, \eqref{eq:rank-condition} is
	equivalent to the following statement: for every $i$, there exist two linearly independent subfamilies that cover the whole family of polynomials and intersect only at $p_i$.
	This can be derived from the finite vector form of Edmonds' matroid partition theorem \cite{Edmonds1965}, with a convenient sufficient condition.

	\begin{corollary}\label{cor:half-independent}
		Let $D\ge1$ and $\ell\ge3$, and let $s=s(D,\ell)$ be as in
		Theorem~\ref{thm:phase-reduction}. If $p_1,\ldots,p_\ell\in\Z[t]$
		have degree at most $D$, zero constant terms, and every
		$\lfloor\ell/2\rfloor+1$ of them are linearly independent over $\Q$,
		then $\W_s(T_1),\ldots,\W_s(T_\ell)$ are characteristic for both
		\eqref{eq:general-average} and \eqref{eq:prime-average}.
		In particular, it suffices that every $\ell-1$ polynomials are
		linearly independent.
	\end{corollary}
	
	\begin{proof}
		For each $i$, divide $[\ell]\setminus\{i\}$ into two sets whose sizes differ by at most one, and add $i$ to these two sets. The size of each resulting set is at most $\lfloor\ell/2\rfloor+1$, so the corresponding subfamilies are linearly independent. Applying Lemma~\ref{lem:rank-cover} and
		Theorem~\ref{thm:phase-reduction} yields the conclusion. The last assertion follows from
		$\lfloor\ell/2\rfloor+1\le\ell-1$.
	\end{proof}
	
	\subsection{Optimal orders for three iterates}
	
	For $t,t^2,t^2+t$, by Theorem~\ref{thm:phase-reduction} one can reduce the average to phase
	factors of some finite order. Subsequently, recurrence relations and polynomial identities give the exact order.
	
	\begin{thmx}\label{thm:main}
		Let $T_1,T_2,T_3$ be invertible commuting measure-preserving
		transformations of a standard probability space $(X,\mathcal X,\mu)$.
		Then the factors $\W_4(T_1),\W_2(T_2),\W_2(T_3)$ are characteristic
		for both $A_N$ and $A_{N,\mathbb P}$.
		More precisely, for any $f_1,f_2,f_3\in L^\infty(\mu)$, if
		\[
		\E(f_1\mid\W_4(T_1))=0,
		\quad\text{or}\quad
		\E(f_2\mid\W_2(T_2))=0,
		\quad\text{or}\quad
		\E(f_3\mid\W_2(T_3))=0,
		\]
		then both averages tend to zero in $L^2(\mu)$.
		For both averaging schemes
		if the conclusion is required to hold for all such systems, then none of the three orders can be reduced even if the constraint
		 $T_2=T_3$ is imposed.
	\end{thmx}
	
	\begin{corollary}\label{cor:orders}
		Under the same system assumptions, the factors
		$\HK_4(T_1),\HK_2(T_2),\HK_2(T_3)$ are characteristic for both
		$A_N$ and $A_{N,\mathbb P}$. In particular,
		\[
		\|f_1\|_{U^5(T_1)}=0
		\quad\text{or}\quad
		\|f_2\|_{U^3(T_2)}=0
		\quad\text{or}\quad
		\|f_3\|_{U^3(T_3)}=0
		\]
		implies that both averages tend to zero in $L^2(\mu)$.
		For both averaging schemes
		if the conclusion is required to hold for all such systems, then none of the three orders can be reduced—even if the constraint
		$T_2=T_3$ is imposed.
	\end{corollary}
	
	The upper bound of this corollary is given by
	$\W_d(T)\subseteq\HK_d(T)$.
	For the lower bounds, Proposition~\ref{prop:sharpness} gives
	commuting affine transformations on $\T^5$, with $T_2=T_3$ along with functions of modulus one that satisfy
	\[
	\|f_1\|_{U^4(T_1)}=0,\qquad
	\|f_2\|_{U^2(T_2)}=\|f_3\|_{U^2(T_3)}=0,
	\]
	Furthermore, the summand in \eqref{eq:average} remains the same function of modulus one for every $n$.
	Consequently, even without imposing restrictions on the factors associated with the other two coordinates, none of the three factor orders can be reduced.
	Within the hierarchies of Host–Kra factors and phase factors, these orders cannot be lowered if the conclusion is required to hold for all commuting systems.
	In particular, the minimal common Host–Kra order is four.
	
\subsection{Outline of the proofs}

For the polynomial average of integer parameters,
by \cite[Theorem~2.8]{FrantzikinakisKuca2025}, each $f_i$ can be replaced by $\E(f_i\mid\HK_s(T_i))$ without changing the $L^2$ limit of the average. Applying \cite[Theorem~1.2]{JamneshanMachado2026}, we approximate these functions in $L^2$ on  $s$-step nilsystem bundles. Since the starting points $\pi_1(x),\ldots,\pi_\ell(x)$
all depend on the same $x$, we may study the polynomial orbits originating from the point $(\pi_1(x),\ldots,\pi_\ell(x))$.

By Leibman’s theorem \cite{Leibman2005}, the limiting orbit measure is a finite average of Haar measures on compact orbits $Hz_1, \ldots, Hz_r$,
where $H \le \prod_i G_i^\circ$ is a connected subgroup.
For each $i$, the rank condition guarantees the existence of two linearly independent subfamilies of polynomials, both of which contain $p_i$, and their union contains the whole family of polynomials.
Using \cite[Theorem~2.10]{FrantzikinakisKuca2025}, we show that the projection of $H$ onto the two corresponding products of groups are both surjective. By calculating the commutator for each coordinate, we obtain
$\prod_i[G_i^\circ,G_i^\circ]\subseteq H$.
Therefore, we can take the commutator quotient for each coordinate,
and then use the torus characters to obtain the polynomial phase for suitable positive integer powers of each transformation.

For the family of polynomials $t, t^2, t^2+t$, we use Weyl’s theorem, recurrence properties, and interpolation to show that 
the degrees of the three phase polynomials
modulo $\Q$ are at most $4$, $2$, and $2$, respectively, at almost every
point where the average along a fixed residue class has a nonzero limit.
By Lemma~\ref{lem:phaseprojection}, these degree bounds modulo $\Q$ give
the characteristic factors $\W_4(T_1)$, $\W_2(T_2)$, and $\W_2(T_3)$.
Finally, using an example involving commuting affine transformations,
it is shown that none of these three orders can be reduced.

For polynomial averages over prime numbers, using 
\cite[Proposition~3.6]{FrantzikinakisHostKra2013}, the proof can be reduced to the case of integer averages over arithmetic progressions.
Since the restriction to an arithmetic sequence does not change the rank of any family of polynomials, in general we still only need the order $s(D,\ell)$. For $t,t^2,t^2+t$, by constructing a finite cyclic extension, we can choose terms in each residue class, thus obtaining the same upper bound of order $4,2,2$.
Furthermore, since the term being averaged in the aforementioned affine transformation example is independent of $n$, this example also proves the optimality of the three orders in the case of prime averaging.

\subsection{Organization}
In Section~\ref{sec:preliminaries}, we recall some definitions and related
results which will be used later.
In Section~\ref{sec:affine}, we prove  the
general integer theorem.
Sections~\ref{sec:interpolation} 
proves the upper bounds $4,2,2$, and
Section~\ref{sec:sharpness}  proves their optimality.
In Section~\ref{sec:primes} we treat prime averages and the
intersection limits.
In Section~\ref{sec:applications} we give further bounds and
examples, including the obstruction at the rank boundary.
	
\subsection*{AI-use disclosure}
AI-assisted tools were used to help extend the results for polynomial averages over the integers to averages over the primes. The paper was written by the author. The author has verified all mathematical
arguments and takes full responsibility for the
correctness of all proofs in the paper.

\subsection*{Acknowledgments}
The author is grateful to his advisor, Professor  Siming Tu,
for introducing him to the subject and for providing
continued guidance, support, and encouragement. 

	\section{Preliminaries}\label{sec:preliminaries}
	
		We write $\N=\{1,2,\ldots\}$, $\T=\R/\Z$, and
	$\mathbb S^1=\{z\in\C:|z|=1\}$, and put
	$e(t)=\exp(2\pi it)$ and $\roots=\{e(t):t\in\Q\}$. For a probability measure $\nu$, let $L^p(\mathcal A, \nu)$ denote the space of all $\mathcal A$-measurable functions in $L^p(\nu)$. When the measure is clear from the context, we use the shorthand notation $\|f\|_p$ and $L^p(\mathcal A)$.
	We denote the conditional expectation by $\E(f\mid\mathcal A)$; in $L^2(\nu)$, this corresponds to the orthogonal projection onto $L^2(\mathcal A, \nu)$.
	The inner product is linear in its first argument. Let $L^\infty(\nu)$ be the space of essentially bounded measurable functions, with the norm $\|f\|_\infty = \operatorname*{ess\,sup}_{\nu}|f|$. For a locally compact space $M$, let $\mathcal B(M)$ denote its Borel $\sigma$-algebra, $C(M)$ the space of complex-valued continuous functions on $M$, and $C_c(M)$ the subspace consisting of functions with compact support. For a Borel probability measure $\mu$ on $M$, define
	$\supp\mu
	=\{x\in M:\mu(U)>0\text{ for every open neighbourhood }U\text{ of }x\}$. We say that $\mu$ has full support if $\supp\mu=M$.
	
	\subsection{Measure-preserving systems and factors}
	
	A measure-preserving system  is a quadruple $(X,\mathcal X,\mu,T)$ in which $(X,\mathcal X,\mu)$ is a standard probability space and $T:X\to X$
	is an invertible measure-preserving transformation. We write $(X,\mu,T)$ when the sigma-algebra is understood.
	We use the convention $Tf=f\circ T$.
	All sets, functions, and maps are measurable.
	In the standard Borel model, functions and maps are Borel-measurable.
	Equalities, inequalities, and invariance are understood modulo
	null sets unless explicitly stated pointwise.
	
	A factor of $(X,\mathcal X,\mu,T)$ is a completed
	$T$-invariant sub-sigma-algebra of $\mathcal X$.
	A factor map $\pi:(X,\mu,T)\to(Y,\nu,S)$ satisfies
	$\pi_*\mu=\nu$ and $\pi\circ T=S\circ\pi$.
	The corresponding factor sigma-algebra on $X$ is the completion
	of $\pi^{-1}\mathcal Y$, where $\mathcal Y$ is the sigma-algebra
	of $Y$. We identify functions on $Y$ with their pullbacks to $X$.
	For a family $\mathcal F$ of measurable sets or functions,
	$\sigma_\nu(\mathcal F)$ denotes the $\nu$-completion of the
	sigma-algebra generated by $\mathcal F$.
	All joins of sigma-algebras are understood to be completed.
	
	Let $\I(T)$ denote the sigma-algebra of sets $A\in\mathcal X$
	such that $\mu(T^{-1}A\mathbin{\triangle}A)=0$.
	The transformation $T$ is ergodic if every $A\in\I(T)$
	has measure zero or one.
	It is totally ergodic if $T^q$ is ergodic for every integer $q\ge1$.
	The rational Kronecker factor of $T$ is
	\[
	\Krat(T)
	=\bigvee_{m\ge1}\I(T^m)
	=\bigvee_{j\ge1}\I(T^{j!}).
	\]

	\subsection{Host--Kra seminorms and phase factors}
	Let $(X,\mathcal X,\mu,T)$ be a	measure-preserving system.
	For  $f\in L^\infty(\mu)$, write
	\(\Delta_Tf=(Tf)\overline f\), \(\Delta_T^0f=f\), and
	\(\Delta_T^{k+1}f=\Delta_T(\Delta_T^kf)\).
	The Host--Kra seminorms satisfy
	\[
	\norm{f}_{U^1(T)}=\norm{\E(f\mid\I(T))}_2,
	\qquad
	\norm{f}_{U^{r+1}(T)}^{2^{r+1}}
	=\lim_{H\to\infty}\av{H}{h}
	\norm{\Delta_{T^h}f}_{U^r(T)}^{2^r}.
	\]
	We use their characterization of \(\HK_d(T)\); see
	\cite{HostKra2005,HostKra2018} and
	\cite[Section~3.2.1]{FrantzikinakisKuca2025} for details.
	
	The mean ergodic theorem gives
	\[
	\lim_{H\to\infty}\av{H}{h}\int\Delta_{T^h}f\,d\mu
	=\int|\E(f\mid\I(T))|^2\,d\mu
	=\norm{f}_{U^1(T)}^2.
	\]
	By induction,
	\begin{equation}\label{eq:iterated-correlation}
		\norm{f}_{U^r(T)}^{2^r}
		=\lim_{H_1\to\infty}\av{H_1}{h_1}\cdots
		\lim_{H_r\to\infty}\av{H_r}{h_r}
		\int\Delta_{T^{h_r}}\cdots\Delta_{T^{h_1}}f\,d\mu.
	\end{equation}
	For \(q\ge1\) and \(r\ge2\), 
	\[
	\norm{f}_{U^r(T)}
	\le\norm{f}_{U^r(T^q)}
	\le q^{r/2^r}\norm{f}_{U^r(T)};
	\]
	see \cite[Lemma~3.1]{FrantzikinakisKuca2025}. Moreover, for any $d\ge1$,
	$\HK_d(T^q)=\HK_d(T)$.
	
	\begin{definition}\label{def:phase-factors}
		Let $(X,\mathcal X,\mu,T)$ be a measure-preserving system and let $d\ge1$
		be an integer. Define
		\[
		\mathcal P_d(T)=\left\{\phi\in L^\infty(\mu):|\phi|=1,\quad
		\Delta_{T^q}^{d+1}\phi=1\text{ for some integer }q\ge1\right\}.
		\]	
		The \emph{phase factor for finite powers} of order $d$ is
		\[
		\W_d(T)=\sigma_\mu(\mathcal P_d(T)).
		\]
		Denote $\mathcal A_d(T)=\operatorname{span}_{\C}\mathcal P_d(T)$.
	\end{definition}

	We next recall Newton's formula for functions
	satisfying $\Delta_S^L\phi=1$.
	For $a\in\R$, write $[a]=a+\Q\in\R/\Q$.
	For $P(t)=\sum_k a_kt^k\in\R[t]$, put
	$[P](t)=\sum_k[a_k]t^k$.
	Here $\R/\Q$ is regarded as a vector space over
	$\Q$, and we set $\deg0=-\infty$.
	This notation will also be used in
	Section~\ref{sec:interpolation}.

\begin{lemma}\label{lem:polynomial-phases}
	Let $(X,\mathcal X,\mu,T)$ be a measure-preserving system, let $L\ge1$
	be an integer, and let
	$\phi\in L^\infty(\mu)$ satisfy
	$|\phi|=1$ and $\Delta_T^L\phi=1$.
	Then there are a measurable $T$-invariant set $X_{T,\phi}$ with $\mu(X_{T,\phi})=1$
	and real polynomials $P_{T,\phi,x}$ of degree less than $L$, with
	coefficients measurable in $x$, such that the following hold.
	\begin{enumerate}
		\item[(i)] For every $x\in X_{T,\phi}$ and $n\in\Z$,
		\begin{equation*}
			\phi(T^nx)
			=\prod_{k=0}^{L-1}\bigl(\Delta_T^k\phi(x)\bigr)^{\binom nk}
			=e(P_{T,\phi,x}(n)).
		\end{equation*}
		
		\item[(ii)] For every $x\in X_{T,\phi}$, $h\in\Z$, and integer $d\ge0$, we have $[P_{T,\phi,T^hx}](t)=[P_{T,\phi,x}](t+h)$ and
		\begin{equation*}
			\begin{aligned}
				\deg[P_{T,\phi,x}]\le d
				\quad\Longleftrightarrow\quad
				\Delta_T^k\phi(x)\in\roots\text{ for every }d<k<L.
			\end{aligned}
		\end{equation*}
		If $d\ge L-1$, both sides of the equivalence hold automatically.
	\end{enumerate}
\end{lemma}

\begin{proof}
	Set
	$\psi_k=\Delta_T^k\phi$ for $0\le k\le L$, and put
	\[
	X_{T,\phi}=\bigcap_{n\in\Z}T^{-n}
	\{x:|\phi(x)|=1,\ \psi_L(x)=1\}.
	\]
	Then $X_{T,\phi}$ is $T$-invariant and $\mu(X_{T,\phi})=1$.
	For every $x\in X_{T,\phi}$ and $0\le k<L$, we have
	\[
	|\psi_k(x)|=1,\qquad
	\psi_k(Tx)=\psi_k(x)\psi_{k+1}(x),\qquad \psi_L(x)=1.
	\]
	Choose measurable $a_k(x)\in[0,1)$ with $e(a_k(x))=\psi_k(x)$ on
	$X_{T,\phi}$, and set $a_k=0$ elsewhere. Define
	\[
	P_{T,\phi,x}(t)=\sum_{k=0}^{L-1}a_k(x)\binom tk,
	\qquad
	\binom t0=1,\quad
	\binom tk=\frac{t(t-1)\cdots(t-k+1)}{k!}\quad(k\ge1).
	\]
	These polynomials have degree less than $L$ and measurable coefficients.
	The formula below also follows from
	\cite[(2.1)--(2.2)]{HaaseMoriakov2018}.
	Here, for every $n\in\Z$, by Pascal's identity we have
	\[
	\begin{aligned}
		\prod_{k=0}^{L-1}\psi_k(Tx)^{\binom nk}
		=\prod_{k=0}^{L-1}
		\bigl(\psi_k(x)\psi_{k+1}(x)\bigr)^{\binom nk}
		=\prod_{k=0}^{L-1}\psi_k(x)^{\binom{n+1}{k}}.
	\end{aligned}
	\]
	The product at $n=0$ is $\phi(x)$. Since $T$ is invertible, by induction one can prove (i) for all $n\in\Z$.
	
	For $x\in X_{T,\phi}$ and $h,n\in\Z$, part (i) gives
	\[
	e(P_{T,\phi,T^hx}(n))
	=\phi(T^{n+h}x)
	=e(P_{T,\phi,x}(n+h)).
	\]
	Since the degree of the polynomial $P_{T,\phi,T^hx}(t)-P_{T,\phi,x}(t+h)$ is  less than $L$ and the polynomial takes
	integer values at $0,1,\ldots,L-1$, Lagrange interpolation
	shows that all its coefficients are rational.
	This proves $[P_{T,\phi,T^hx}](t)=[P_{T,\phi,x}](t+h)$.
	For the remaining part, use
	\[
	[P_{T,\phi,x}](t)=\sum_{k=0}^{L-1}[a_k(x)]\binom tk,
	\qquad
	[a_k(x)]=0\quad\Longleftrightarrow\quad\psi_k(x)\in\roots.
	\]
	If $[a_k(x)]=0$ for all $d<k<L$, then
	$\deg[P_{T,\phi,x}]\le d$.
	
	Conversely, assume that $[a_k(x)]\ne0$ for some $k>d$.
	Take the largest such $k$.
	Since $\binom tk$ has leading coefficient $1/k!$,
	the coefficient of $t^k$ in $[P_{T,\phi,x}]$ is
	$[a_k(x)]/k!\ne0$. Hence
	\[
	\deg[P_{T,\phi,x}]\ge k>d.
	\]
\end{proof}
	\begin{remark}
For a measure-preserving system $(X,\mathcal X,\mu,T)$ , define the Abramov factor
\[
\operatorname{Abr}_{<d+1}(T)
=\sigma_\mu\{\phi\in L^\infty(\mu):
|\phi|=1,\ \Delta_T^{d+1}\phi=1\}.
\]
This agrees with the mixed-difference definition in
\cite[Definition~1.13]{BergelsonTaoZiegler2010}.
Indeed, Lemma~\ref{lem:polynomial-phases} shows that
$\Delta_T^{d+1}\phi=1$ implies
$\Delta_{T^{n_1}}\cdots\Delta_{T^{n_{d+1}}}\phi=1$
for all $n_1,\ldots,n_{d+1}\in\Z$.
The converse follows by taking $n_1=\cdots=n_{d+1}=1$.

The phase factor $\W_d(T)$ allows the same difference identity
to hold for a positive power of $T$, with the power depending
on the generating function. Thus
\[
\W_d(T)=\bigvee_{q\ge1}\operatorname{Abr}_{<d+1}(T^q).
\]
In particular, $\operatorname{Abr}_{<d+1}(T)\subseteq\W_d(T)$.
	\end{remark}
	
	The following lemma gives the properties of
	$\W_d(T)$ needed in the proofs below.

	\begin{lemma}\label{lem:phasefactors}
		Let $(X,\mathcal X, \mu,T)$ be a measure-preserving system and let $d\ge1$
		be an integer. Then the following hold.
		\begin{enumerate}
			\item[(i)]  If $S$ is
			an invertible measure-preserving transformation of $(X,\mathcal X,\mu)$
			commuting with $T$,
			then
			$S^{-1}\W_d(T)=\W_d(T)$.
			In particular, the sigma-algebra $\W_d(T)$ is $T$-invariant.
			\item[(ii)] For every $q\ge1$,
			\[
			\W_d(T^q)=\W_d(T),\qquad
			\I(T^q)\subseteq\W_d(T)\subseteq\HK_d(T).
			\]
			
			\item[(iii)] For every $f\in L^2(\W_d(T))$ and $\varepsilon>0$,
			there exists $U\in\mathcal A_d(T)$ such that
			\[
			\|f-U\|_2<\varepsilon.
			\]
		\end{enumerate}
	\end{lemma}
	
	\begin{proof}
		Fix $\phi\in\mathcal P_d(T)$. Choose $m\ge1$ such that
		$\Delta_{T^m}^{d+1}\phi=1$, and put $R=T^m$.
		By Lemma~\ref{lem:polynomial-phases}(i), for almost every $x$ and every $n\in\Z$, we have
		\[
		\phi(R^n x)=e(P_{R,\phi,x}(n)),
		\qquad \deg P_{R,\phi,x}\le d.
		\]
		For every $h\ge1$,
		$\Delta_{R^h}\phi(R^n x)
		=e(P_{R,\phi,x}(n+h)-P_{R,\phi,x}(n))$.
		Since $\deg P_{R,\phi,x}\le d$, taking $d+1$ differences gives
		$\Delta_{R^{h_{d+1}}}\cdots\Delta_{R^{h_1}}\phi=1$
		for all $h_1,\ldots,h_{d+1}\ge1$.
		In particular,
		$\Delta_{R^q}^{d+1}\phi=1$ for every $q\ge1$.
		
		For (i), suppose that $S$ commutes with $T$.
		For $\phi\in\mathcal P_d(T)$, choose $m\ge1$ such that
		$\Delta_{T^m}^{d+1}\phi=1$. Then
		$\Delta_{T^m}^{d+1}(\phi\circ S^{\pm1})
		=(\Delta_{T^m}^{d+1}\phi)\circ S^{\pm1}=1$.
		Thus $\phi\circ S$ and $\phi\circ S^{-1}$ belong to
		$\mathcal P_d(T)$. Since $\W_d(T)=\sigma(\mathcal P_d(T))$,
		$S^{-1}\W_d(T)=\W_d(T)$.
		Taking $S=T$ shows that $\W_d(T)$ is $T$-invariant.
		
		For (ii), by definition, we have
		$\mathcal P_d(T^q)\subseteq\mathcal P_d(T)$.
		Conversely, if $\Delta_{T^m}^{d+1}\phi=1$, the argument above gives
		$\Delta_{T^{mq}}^{d+1}\phi=1$.
		Therefore
		$\mathcal P_d(T^q)=\mathcal P_d(T)$ and
		$\W_d(T^q)=\W_d(T)$.
		If $B\in\I(T^q)$, put $\phi=1-2\,1_B$. Then
		$|\phi|=1$ and $\Delta_{T^q}\phi=1$.
		Thus $\phi\in\mathcal P_d(T)$, and
		$1_B=\frac{1-\phi}{2}\in L^\infty(\W_d(T))$.
		This proves $\I(T^q)\subseteq\W_d(T)$.
		
		It remains to prove $\W_d(T)\subseteq\HK_d(T)$.
		Fix $\phi\in\mathcal P_d(T)$ and choose $m\ge1$ with
		$\Delta_{T^m}^{d+1}\phi=1$. Put $R=T^m$ and
		$g=\E(\phi\mid\HK_d(R))$.
		Since $\Delta_R^{d+1}\phi=1$, the preceding argument gives
		$\Delta_{R^{h_{d+1}}}\cdots\Delta_{R^{h_1}}\phi=1$
		for all $h_1,\ldots,h_{d+1}\ge1$.
		By \eqref{eq:iterated-correlation},
		$\|\phi\|_{U^{d+1}(R)}=1$.
		Since $g=\E(\phi\mid\HK_d(R))$, we have
		$\|\phi-g\|_{U^{d+1}(R)}=0$. Thus
		\[
		1=\|\phi\|_{U^{d+1}(R)}
		\le\|g\|_{U^{d+1}(R)}+\|\phi-g\|_{U^{d+1}(R)}
		=\|g\|_{U^{d+1}(R)}
		\le\|\phi\|_{U^{d+1}(R)}=1.
		\]
		Hence $\|g\|_{U^{d+1}(R)}=1$.
		
		Since $|\phi|=1$, conditional expectation gives $|g|\le1$.
		For all $h_1,\ldots,h_{d+1}\ge1$,
		\[
		\begin{aligned}
			\left|\Delta_{R^{h_{d+1}}}\cdots
			\Delta_{R^{h_1}}g(x)\right|
			=
			\prod_{\omega\in\{0,1\}^{d+1}}
			\left|g\left(R^{\omega_1h_1+\cdots+\omega_{d+1}h_{d+1}}x\right)\right|
			\le |g(x)|.
		\end{aligned}
		\]
		Therefore, for every $H_1,\ldots,H_{d+1}\ge1$,
		\[
		\begin{aligned}
			&\quad\left|
			\frac1{H_1\cdots H_{d+1}}
			\sum_{h_1=1}^{H_1}\cdots\sum_{h_{d+1}=1}^{H_{d+1}}
			\int_X\Delta_{R^{h_{d+1}}}\cdots
			\Delta_{R^{h_1}}g\,d\mu
			\right|\\
			&\le
			\frac1{H_1\cdots H_{d+1}}
			\sum_{h_1=1}^{H_1}\cdots\sum_{h_{d+1}=1}^{H_{d+1}}
			\int_X
			\left|\Delta_{R^{h_{d+1}}}\cdots
			\Delta_{R^{h_1}}g\right|\,d\mu\\
			&\le\int_X|g|\,d\mu.
		\end{aligned}
		\]
		Taking the successive limits in
		\eqref{eq:iterated-correlation}, we obtain
		$1=\|g\|_{U^{d+1}(R)}^{2^{d+1}}
		\le\int_X|g|\,d\mu\le1$.
		Thus $|g|=1$. Since $g$ is the orthogonal
		projection of $\phi$ onto $L^2(\HK_d(R))$,
		$\int_X|\phi-g|^2\,d\mu
		=\int_X|\phi|^2\,d\mu-\int_X|g|^2\,d\mu
		=0$.
		Thus $\phi\in L^2(\HK_d(R))$. Since
		\[
		\|f\|_{U^{d+1}(T)}
		\le \|f\|_{U^{d+1}(T^m)}
		\le m^{(d+1)/2^{d+1}}\|f\|_{U^{d+1}(T)},
		\]
		then
		$\HK_d(T^m)=\HK_d(T)$.
		Every $\phi\in\mathcal P_d(T)$
		is therefore $\HK_d(T)$-measurable, proving
		$\W_d(T)\subseteq\HK_d(T)$.
		
		For (iii), let $\phi,\chi\in\mathcal P_d(T)$.
		Choose $m_1,m_2\ge1$ such that
		$\Delta_{T^{m_1}}^{d+1}\phi=1$ and
		$\Delta_{T^{m_2}}^{d+1}\chi=1$.
		Taking $q=m_1m_2$, we have
		$\Delta_{T^q}^{d+1}\phi
		=\Delta_{T^q}^{d+1}\chi=1$.
		Hence
		$\Delta_{T^q}^{d+1}(\phi\chi)=1$ and
		$\Delta_{T^q}^{d+1}\overline\phi=1$.
		Thus $\phi\chi,\overline\phi\in\mathcal P_d(T)$.
		Since $1\in\mathcal P_d(T)$,
		$\mathcal A_d(T)$ is an algebra containing the constants
		and closed under complex conjugation.
		Moreover,
		$\sigma(\mathcal A_d(T))=\W_d(T)$.
		By the monotone class theorem, every bounded
		$\W_d(T)$-measurable function can be approximated in
		$L^2$ by functions in $\mathcal A_d(T)$.
		Since bounded functions are dense in $L^2(\W_d(T))$,
		the proof is complete.
	\end{proof}

	We need the following $L^2$ approximation lemma.

	\begin{lemma}\label{lem:l2extension}
		Let $(X,\mathcal X,\mu)$ be a probability space, let $\ell\ge1$,
		and let $B_N:(L^\infty(\mu))^\ell\to L^2(\mu)$, $N\ge1$,
		be complex multilinear maps. Suppose that a constant $C\ge0$ satisfies
		\[
		\|B_N(f_1,\ldots,f_\ell)\|_2
		\le C\|f_i\|_2\prod_{j\ne i}\|f_j\|_\infty
		\]
		for all bounded $f_1,\ldots,f_\ell$, $i\in[\ell]$, and $N\ge1$.
		For each $i$, let $\mathcal V_i\subseteq\mathcal H_i\subseteq L^\infty(\mu)$
		be complex linear spaces, with $\mathcal V_i$ dense in $\mathcal H_i$
		in $L^2$.
		If $\|B_N(v_1,\ldots,v_\ell)\|_2\to0$ whenever
		$v_i\in\mathcal V_i$, then
		$\|B_N(f_1,\ldots,f_\ell)\|_2\to0$ whenever $f_i\in\mathcal H_i$.
	\end{lemma}
	
	\begin{proof}
		We extend the assertion one coordinate at a time.
		Suppose it holds when the first $r-1$ coordinates belong to
		$\mathcal H_i$ and the remaining coordinates belong to $\mathcal V_i$.
		Fix $f_i\in\mathcal H_i$ for $i\le r$ and
		$v_i\in\mathcal V_i$ for $i>r$.
		For every $v_r\in\mathcal V_r$, multilinearity gives
		\[
		\limsup_{N\to\infty}
		\|B_N(f_1,\ldots,f_r,v_{r+1},\ldots,v_\ell)\|_2
		\le C\|f_r-v_r\|_2
		\prod_{i<r}\|f_i\|_\infty\prod_{i>r}\|v_i\|_\infty.
		\]
		Approximate $f_r$ in $L^2$ by elements of $\mathcal V_r$ while
		keeping the other functions fixed. The right-hand side tends to zero.
		Starting at $r=1$ and continuing to $r=\ell$ proves the result.
	\end{proof}
	
	\subsection{Nilmanifolds and nilsystems}
	For a group $G$, write $1_G$ for its identity and
	$[a,b]=aba^{-1}b^{-1}$ for the commutator of $a,b\in G$.
	For subgroups $A,B\le G$, let $[A,B]$ be the subgroup generated
	by their commutators. The lower central series is defined by
	$\gamma_1(G)=G$ and $\gamma_{r+1}(G)=[G,\gamma_r(G)]$ for $r\ge1$.
	For an integer $s\ge1$, we call $G$ $s$-step nilpotent if $\gamma_{s+1}(G)=\{1_G\}$. For a Lie group $G$,
	$G^\circ$ denotes the connected component of $1_G$.
	
	Let $G$ be an $s$-step nilpotent Lie group and let $\Gamma\le G$
	be a discrete cocompact subgroup. Then $M=G/\Gamma$
	is called an $s$-step nilmanifold. We write $gx$ for the left action
	$g(h\Gamma)=gh\Gamma$ and $m_M$ for the unique $G$-invariant
	Borel probability measure on $M$. For $\tau\in G$, define
	$R_\tau:x\mapsto\tau x$ $(x\in M)$ and the system
	$(M,m_M,R_\tau)$ is called an $s$-step nilsystem. 
	
	For a nilmanifold $M$,
	$\operatorname{Comp}(M)$ denotes the finite set
	of connected components of $M$.
	Any $R\in\operatorname{Comp}(M)$ is clopen and can
	be written as $R=G^\circ g_R\Gamma/\Gamma$ for some
	$g_R\in G$. Let $m_R=(m_M|_R)/m_M(R)$ be its normalized Haar measure.
	For a closed subgroup $H\le G$ with a compact
	orbit $Hx$, we write $m_{Hx}$ for the $H$-invariant probability
	measure on that orbit. A continuous character of a torus $K$ is
	a continuous homomorphism $\chi:K\to\mathbb S^1$.
	
	Let $C=[G^\circ,G^\circ]$. The \emph{commutator quotient map} is
	$\rho:M\to Z$, where $Z=\{Cx:x\in M\}$ with the quotient topology,
	$Cx=\{cx:c\in C\}$, and $\rho(x)=Cx$.
	For any connected component $R=G^\circ g_R\Gamma/\Gamma$ of $M$, let
	$\Lambda_R=G^\circ\cap g_R\Gamma g_R^{-1}$.
	The quotient $Z_R=\rho(R)$ is identified with a torus by
	\[
	Z_R\simeq G^\circ/(C\Lambda_R)\simeq\T^{d_R},\qquad
	\rho(hg_R\Gamma)=h(C\Lambda_R);
	\]
	see for example \cite[Section~1.10]{Leibman2005}.
	Hence $Z=\bigsqcup_{R\in\operatorname{Comp}(M)}Z_R$ and
	$\rho^{-1}(\rho(x))=Cx$.

\section{Reduction to phase factors}\label{sec:affine}

In this section, we prove the integer assertion of Theorem~\ref{thm:phase-reduction}. We first express the rank condition in terms
of two linearly independent subfamilies, then use nilsystem
models of the individual Host--Kra factors.
	
	\subsection{The rank condition and independent subfamilies}
	
	For a family of polynomials $\mathbf{p} = (p_1, \ldots, p_\ell)$ with a finite index set,
	let $r_{\mathbf{p}}(J) = \dim_{\mathbb{Q}} \operatorname{span}_{\mathbb{Q}} \{p_j : j \in J\}$.
	The following lemma provides an equivalent form of the rank condition.

	\begin{lemma}\label{lem:rank-cover}
		Let $\ell\ge1$ and $p_1,\ldots,p_\ell\in\Q[t]$.
		The following are equivalent.
		\begin{enumerate}
			\item[(i)] For every nonempty $J\subseteq[\ell]$,
			\[
			|J|<2r_{\mathbf p}(J).
			\]
			
			\item[(ii)] For every $i\in[\ell]$, there exist
			$J_i,K_i\subseteq[\ell]$ such that
			\[
			i\in J_i\cap K_i,\qquad
			J_i\cup K_i=[\ell],
			\]
			and both $(p_j)_{j\in J_i}$ and
			$(p_j)_{j\in K_i}$ are linearly independent
			over $\Q$.
		\end{enumerate}
	\end{lemma}
	
	\begin{proof}
		Assume (ii) holds. Fix a nonempty
		$J\subseteq[\ell]$ and choose $i\in J$.
		Since $J_i\cup K_i=[\ell]$,
		$i\in J_i\cap K_i$, and both $(p_j)_{j\in J_i}$ and
		$(p_j)_{j\in K_i}$ are linearly independent
		over $\Q$ we have
		\[
		\begin{aligned}
			|J|
			&\le |J\cap J_i|+|J\cap K_i|-1\\
			&\le 2r_{\mathbf p}(J)-1
			<2r_{\mathbf p}(J).
		\end{aligned}
		\]
		This proves (i).
		
		Conversely, suppose (i) holds and fix
		$i\in[\ell]$. Take $J=\{i\}$ with
		$p_i\ne0$. Consider the quotient vector space
		\[
		V=\operatorname{span}_{\Q}
		\{p_1,\ldots,p_\ell\}/\Q p_i,
		\]
		and  let
		$\overline p_j=p_j+\Q p_i\in V$ for $j\ne i$.
		Then for every $A\subseteq[\ell]\setminus\{i\}$,
		\[
		\dim_{\Q}\operatorname{span}_{\Q}
		\{\overline p_j:j\in A\}
		=r_{\mathbf p}(A\cup\{i\})-1.
		\]
		Applying (i) to $A\cup\{i\}$ gives
		$|A|+1<2r_{\mathbf p}(A\cup\{i\})$.
		Since both sides are integers,
		\[
		|A|\le2\bigl(r_{\mathbf p}(A\cup\{i\})-1\bigr)
		=2\dim_{\Q}\operatorname{span}_{\Q}
		\{\overline p_j:j\in A\}.
		\]
		By Edmonds's partition theorem
		\cite[Theorem~1]{Edmonds1965}, the vectors
		$(\overline p_j)_{j\ne i}$ can therefore be
		partitioned into two linearly independent
		subfamilies. Thus
		\[
		[\ell]\setminus\{i\}=A_i\sqcup B_i,
		\]
		where $(\overline p_j)_{j\in A_i}$ and
		$(\overline p_j)_{j\in B_i}$ are both
		linearly independent.
		
		Put $J_i=A_i\cup\{i\}$ and
		$K_i=B_i\cup\{i\}$. Then
		$J_i\cup K_i=[\ell]$,
		$J_i\cap K_i=\{i\}$.
		To prove that $(p_j)_{j\in J_i}$ is linearly
		independent, suppose
		\[
		c_ip_i+\sum_{j\in A_i}c_jp_j=0,
		\qquad c_i,c_j\in\Q.
		\]
		we have
		\[
		\sum_{j\in A_i}c_j\overline p_j
		=\sum_{j\in A_i}c_jp_j+\Q p_i
		=-c_ip_i+\Q p_i=0
		\quad\text{in }V.
		\]
		As $(\overline p_j)_{j\in A_i}$ is linearly
		independent, $c_j=0$ for all $j\in A_i$.
		Since $p_i\ne0$, we also have $c_i=0$.
		Thus $(p_j)_{j\in J_i}$ is independent.
		The same argument applies to $K_i$,
		proving (ii).
	\end{proof}

	\subsection{Nilsystem bundles and commutator quotients}
	
	\begin{definition}[Nilsystem bundles]\label{def:nilsystem-bundle}
		Following \cite[Definition~2.2]{JamneshanMachado2026},
		A \emph{nilsystem bundle} is a measure-preserving system $(Y,\nu,S)$
		of the form
		\[
		\begin{gathered}
			Y=\bigsqcup_{a\in I}\Omega_a\times M_a,\qquad
			M_a=G_a/\Gamma_a,\qquad I\subseteq\N,\\
			\nu=\sum_{a\in I}\beta_a\otimes m_a,\qquad
			S(\omega,y)=(\omega,\tau_a(\omega)y)
			\quad\text{on }\Omega_a\times M_a.
		\end{gathered}
		\]
		Every $\Omega_a$ is a standard Borel space with a finite measure
		$\beta_a$, and $\sum_{a\in I}\beta_a(\Omega_a)=1$.
		Every $M_a$ is a nilmanifold with Haar measure
		$m_a=m_{M_a}$. The maps $\tau_a:\Omega_a\to G_a$ are measurable,
		and $y\mapsto\tau_a(\omega)y$ is ergodic on $(M_a,m_a)$ for
		$\beta_a$-almost every $\omega$.
		The sigma-algebra $\mathcal Y$ on $Y$ is the $\nu$-completion of
		the disjoint-union product sigma-algebra.
		
		For an integer $s\ge1$, the bundle is called \emph{$s$-step}
		if every $G_a$ has nilpotency class at most $s$.
		Without this qualification, no common bound on these classes is assumed.
	\end{definition}
	
	For several nilsystem bundles $(Y_i,\nu_i,S_i)$, we use an additional
	index $i$, writing $\Omega_{i,a}$, $M_{i,a}$,
	$\beta_{i,a}$, and $\tau_{i,a}$.
	
	\begin{lemma}
		\label{lem:commutator-quotient}
		Let $M=G/\Gamma$ be a nilmanifold with Haar measure $m_M$ and put
		$C=[G^\circ,G^\circ]$.
		Let $Z=\{Cx:x\in M\}$ with the quotient topology, and define
		$\rho:M\to Z$ by $\rho(x)=Cx=\{cx:c\in C\}$.
		Let $m_Z=\rho_*m_M$.
		If $g\in L^2(m_M)$ satisfies
		$\|g(c\,\cdot)-g\|_2=0$ for every $c\in C$, then there exists
		$\widetilde g\in L^2(m_Z)$ such that
		$g=\widetilde g\circ\rho$ $m_M$-almost everywhere.
	\end{lemma}
	
	\begin{proof}
		First assume that $|g|\le K$ almost everywhere and fix a left Haar measure $dh$ on $G^\circ$.
		Take a decreasing neighbourhood basis \((V_n)\) at the identity
		of \(G^\circ\), and \(\eta_n\in C_c(G^\circ)\) with
		\[
		\eta_n\ge0,\qquad \int\eta_n(h)\,dh=1,\qquad
		\supp\eta_n\subseteq V_n,\qquad
		g_n(x)=\int_{G^o}\eta_n(h)g(h^{-1}x)\,dh.
		\]
		For every $F\in C(M)$, by the $G^\circ$-invariance of $m_M$
		\[
		\|g(h^{-1}\,\cdot)-g\|_2
		\le 2\|g-F\|_2+
		\sup_{x\in M}|F(h^{-1}x)-F(x)|.
		\]
		For fixed $F$, compactness gives
		$\limsup_{h\to1_G}\|g(h^{-1}\,\cdot)-g\|_2\le2\|g-F\|_2$.
		Approximating $g$ by continuous $F$ therefore gives
		$\|g(h^{-1}\,\cdot)-g\|_2\to0$ as $h\to1_G$. Hence
		\[
		\|g_n-g\|_2
		\le\sup_{h\in\operatorname{supp}\eta_n}
		\|g(h^{-1}\,\cdot)-g\|_2
		\longrightarrow0.
		\]

		For $a\in G^\circ$, changing variables gives
		\[
		|g_n(ax)-g_n(x)|
		\le K\int|\eta_n(ah)-\eta_n(h)|\,dh
		\longrightarrow0\quad
		\text{when~} a\to1_G,
		\]
		uniformly in $x\in M$.
		Since $Vx$ is open in $M$ for every open
		neighborhood $V$ of $1_G$ in $G^\circ$,
		we have $g_n\in C(M)$.
		For $c\in C=[G^\circ,G^\circ]$, we have
		\[
		\begin{aligned}
			\|g_n(c\,\cdot)-g_n\|_2
			\le\int\eta_n(h)
			\|g((h^{-1}ch)\,\cdot)-g\|_2\,dh
			=0.
		\end{aligned}
		\]
		Since $g_n$ is continuous and $m_M$ has full
		support, the equality
		$\|g_n(c\,\cdot)-g_n\|_2=0$ implies
		$g_n(cx)=g_n(x)$
		for every $c\in C,\ x\in M$.
		Thus the formula
		$\widetilde g_n(Cx)=g_n(x)$
		defines a continuous function on $Z$
		by the definition of the quotient topology.
		Hence $g_n=\widetilde g_n\circ\rho$.
		By $m_Z=\rho_*m_M$,
		$\|\widetilde g_n-\widetilde g_m\|_{L^2(m_Z)}
		=\|g_n-g_m\|_2$.
		Hence $\widetilde g_n\to\widetilde g$ in
		$L^2(m_Z)$ for some $\widetilde g$, and
		\[
		\|g-\widetilde g\circ\rho\|_2
		\le\|g-g_n\|_2+
		\|\widetilde g_n-\widetilde g\|_{L^2(m_Z)}
		\longrightarrow0.
		\]

		For general $g\in L^2(m_M)$, put
		$g_N=g1_{\{|g|\le N\}}$. Then $g_N$ is
		bounded, is $C$-invariant in $L^2$, and
		$\|g-g_N\|_2^2=\int_{\{|g|>N\}}|g|^2\,dm_M\longrightarrow0$, when $N\to\infty$.
		By the bounded case, there exists
		$\widetilde g_N\in L^2(m_Z)$ such that
		$g_N=\widetilde g_N\circ\rho$. Then as $N,L\to\infty$,
		$\|\widetilde g_N-\widetilde g_L\|_{L^2(m_Z)}
		=\|g_N-g_L\|_2\longrightarrow0$.
		Thus $\widetilde g_N$ converges in $L^2(m_Z)$
		to some $\widetilde g$, and
		\[
		\|g-\widetilde g\circ\rho\|_2
		\le\|g-g_N\|_2+
		\|\widetilde g_N-\widetilde g\|_{L^2(m_Z)}
		\longrightarrow0.
		\]

	\end{proof}

	\begin{lemma}
		\label{lem:bundle-phases}\label{lem:commutator-density}
		Let $s\ge1$, and let $(Y,\nu,S)$ be an $s$-step
		nilsystem bundle, where
		$Y=\bigsqcup_{a\in I}\Omega_a\times M_a$,
		$M_a=G_a/\Gamma_a$, $I\subseteq\N$,
		$\nu=\sum_{a\in I}\beta_a\otimes m_a$, and
		$S(\omega,y)=(\omega,\tau_a(\omega)y)$
		on each $\Omega_a\times M_a$.
		For each $a\in I$, put
		$C_a=[G_a^\circ,G_a^\circ]$ and
		$Z_a=\{C_a y:y\in M_a\}$, and define
		$\rho_a:M_a\to Z_a$ by $\rho_a(y)=C_a y$.
		Give $Z_a$ the quotient topology.
		Let
		\[
		\mathcal F=\sigma_\nu\bigl(
		\{B\times\rho_a^{-1}(D):
		a\in I,\ B\subseteq\Omega_a\text{ measurable},
		\ D\in\mathcal B(Z_a)\}\bigr).
		\]
		Define
		\[
		\begin{aligned}
			\mathcal V
			&=\operatorname{span}_{\C}
			\{b(\omega)v(y):
			a\in I,\ b\in L^\infty(\beta_a),\
			v\in C(M_a)\},\\
			\mathcal D
			&=\operatorname{span}_{\C}
			\{b(\omega)(v(cy)-v(y)):
			a\in I,\ b\in L^\infty(\beta_a),\
			v\in C(M_a),\ c\in C_a\},
		\end{aligned}
		\]
		where each function is defined on
		$\Omega_a\times M_a$ and is zero elsewhere.
		
		\begin{enumerate}
			\renewcommand{\labelenumi}{(\roman{enumi})}
			\item $\mathcal F\subseteq\W_s(S)$.
			
			\item $\mathcal V$ is dense in $L^2(\nu)$, and
			$\mathcal D^\perp=L^2(\mathcal F,\nu)$.
			
			\item If $F\in L^2(\nu)$ satisfies
			$\E(F\mid\W_s(S))=0$, then there exist
			$\psi_k\in\mathcal D$ such that
			$\|F-\psi_k\|_{L^2(\nu)}\to0$.
		\end{enumerate}
	\end{lemma}
	
	\begin{proof}
		\noindent\emph{(i).}
		For each $a$, put $q_a=|\operatorname{Comp}(M_a)|!$.
		Then $\tau_a(\omega)^{q_a}R=R$ for every
		$\omega\in\Omega_a$ and $R\in\operatorname{Comp}(M_a)$.
		
		Fix $a$ and a connected component
		$R=G_a^\circ g_R\Gamma_a/\Gamma_a$ of $M_a$.
		Put
		\[
		\Lambda_R=G_a^\circ\cap g_R\Gamma_a g_R^{-1},
		\qquad
		Z_{a,R}=G_a^\circ/(C_a\Lambda_R).
		\]
		Let $\chi$ be a continuous character of $Z_{a,R}$.
		Put $f_0=\chi\circ\rho_a$ on $R$ and
		$\chi_0(h)=\chi(h(C_a\Lambda_R))$ for $h\in G_a^\circ$.
		Since $\rho_a(hg_R\Gamma_a)=h(C_a\Lambda_R)$,
		$f_0(hy)=\chi_0(h)f_0(y)$ for all
		$h\in G_a^\circ,\ y\in R$.

		Fix $\omega\in\Omega_a$ and put
		$u=\tau_a(\omega)^{q_a}$.
		Define $f_{r+1}(y)=f_r(uy)\overline{f_r(y)}$.
		For $h\in G_a^\circ$, define
		$h_0(h)=h$ and $h_{r+1}(h)=[u,h_r(h)]$.
		Then
		$h_r(h)\in G_a^\circ\cap\gamma_{r+1}(G_a)$
		and
		$h_r(uhu^{-1})=uh_r(h)u^{-1}$.
		
		We prove by induction that
		$f_r(hy)=\chi_0(h_r(h))f_r(y)$.
		The case $r=0$ follows from the definition of
		$f_0$. For the induction step,
		\[
		\begin{aligned}
			f_{r+1}(hy)
			&=f_r((uhu^{-1})uy)\overline{f_r(hy)}\\
			&=\chi_0(uh_r(h)u^{-1})
			\overline{\chi_0(h_r(h))}f_{r+1}(y)\\
			&=\chi_0(h_{r+1}(h))f_{r+1}(y).
		\end{aligned}
		\]
		Since $G_a$ has step at most $s$, we have
		$h_s(h)=1_{G_a}$. Hence $f_s(hy)=f_s(y)$
		for every $h\in G_a^\circ$ and $y\in R$.
		Since $R=G_a^\circ g_R\Gamma_a/\Gamma_a$,
		there is $\zeta\in\mathbb S^1$ such that
		$f_s(y)=\zeta$ for all $y\in R$.
		As $uR=R$,
		$f_{s+1}(y)
		=f_s(uy)\overline{f_s(y)}
		=\zeta\overline\zeta=1$ for all $y\in R$.
		
		Define $\Phi$ on $Y$ by
		\[
		\Phi(\omega,y)=
		\begin{cases}
			\chi(\rho_a(y)),
			&(\omega,y)\in\Omega_a\times R,\\
			1,&\text{otherwise}.
		\end{cases}
		\]
		For $(\omega,y)\in\Omega_a\times R$, we have
		$S^{q_a}(\omega,y)=(\omega,uy)$ and $uR=R$.
		Thus $\Delta_{S^{q_a}}^r\Phi(\omega,y)=f_r(y)$
		for every $r\ge0$, and in particular
		$\Delta_{S^{q_a}}^{s+1}\Phi=1$ on
		$\Omega_a\times R$.
		The same holds on $Y\setminus(\Omega_a\times R)$,
		since this set is $S^{q_a}$-invariant and $\Phi=1$ there.
		Hence $\Phi\in L^\infty(\W_s(S))$.
		
		For any measurable $B\subseteq\Omega_a$,
		we have $S^{-q_a}(B\times R)=B\times R$.
		By Lemma~\ref{lem:phasefactors},
		$B\times R\in\I(S^{q_a})\subseteq\W_s(S)$.
		Since $\Phi\in L^\infty(\W_s(S))$, we also have
		$1_{B\times R}\Phi\in L^\infty(\W_s(S))$.
		The continuous characters of the torus
		$Z_{a,R}$ generate its Borel sigma-algebra.
		Consequently,
		$B\times\rho_a^{-1}(D_R)\in\W_s(S)$
		for all measurable $B\subseteq\Omega_a$
		and $D_R\in\mathcal B(Z_{a,R})$.
		
		For every $D\in\mathcal B(Z_a)$, we have
		\[
		B\times\rho_a^{-1}(D)
		=\bigcup_{R\in\operatorname{Comp}(M_a)}
		B\times\rho_a^{-1}(D\cap Z_{a,R})
		\in\W_s(S).
		\]
		Since these sets generate $\mathcal F$ after
		completion, we obtain
		$\mathcal F\subseteq\W_s(S)$.
		
		\smallskip
		\noindent\emph{(ii).}
		For measurable $B\subseteq\Omega_a$ and $E\subseteq M_a$,
		take $v_n\in C(M_a)$ with
		$\|v_n-1_E\|_{L^2(m_a)}\to0$. Then
		\[
		\|1_B(\omega)v_n(y)-1_{B\times E}\|_{L^2(\nu)}^2
		=\beta_a(B)\|v_n-1_E\|_{L^2(m_a)}^2
		\longrightarrow0.
		\]
		Finite linear combinations of rectangle indicators
		on finitely many pieces are dense in $L^2(\nu)$.
		Hence the functions $1_B(\omega)v(y)$,
		with $v\in C(M_a)$, have dense linear span,
		and $\mathcal V$ is dense.
		
		Let $H\in L^2(\mathcal F,\nu)$.
		Since $\rho_a(cy)=\rho_a(y)$,
		$H(\omega,cy)=H(\omega,y)$ almost everywhere
		for each $c\in C_a$. Haar invariance gives
		\[
		\begin{aligned}
			&\langle H,b(\omega)(v(cy)-v(y))\rangle\\
			&\quad=\int_{\Omega_a\times M_a}
			\bigl(H(\omega,c^{-1}y)-H(\omega,y)\bigr)
			\overline{b(\omega)v(y)}\,
			d(\beta_a\otimes m_a)=0.
		\end{aligned}
		\]
		Thus $L^2(\mathcal F,\nu)\subseteq\mathcal D^\perp$.
		
		Since $L^2(\mathcal F,\nu)$ is closed,
		it remains to prove that
		$H\perp\mathcal D$ and $H\perp L^2(\mathcal F,\nu)$
		imply $H=0$.
		Fix $a\in I$ and a measurable $B\subseteq\Omega_a$.
		Put $h_B(y)=\int_B H(\omega,y)\,d\beta_a(\omega)$.
		By Cauchy--Schwarz,
		\[
		\|h_B\|_{L^2(m_a)}^2
		\le\beta_a(B)
		\int_{B\times M_a}|H|^2\,d(\beta_a\otimes m_a)
		<\infty.
		\]
		For $c\in C_a$ and $v\in C(M_a)$,
		\[
		\begin{aligned}
			\langle h_B(c^{-1}\,\cdot)-h_B,v\rangle
			=\langle h_B,v(c\,\cdot)-v\rangle
			=\langle H,1_B(\omega)(v(cy)-v(y))\rangle
			=0.
		\end{aligned}
		\]
		Density of $C(M_a)$ in $L^2(m_a)$ gives
		$h_B(c^{-1}\,\cdot)=h_B$ in $L^2(m_a)$
		for every $c\in C_a$.
		By Lemma~\ref{lem:commutator-quotient},
		$1_B(\omega)h_B(y)\in L^2(\mathcal F,\nu)$.
		Therefore, Fubini's theorem gives
		\[
		\begin{aligned}
			0
			=\langle H,1_B(\omega)h_B(y)\rangle
			&=\int_{M_a}
			\left(\int_B H(\omega,y)\,d\beta_a(\omega)\right)
			\overline{h_B(y)}\,dm_a(y)\\
			&=\|h_B\|_{L^2(m_a)}^2.
		\end{aligned}
		\]
		Hence $h_B=0$, and
		$\langle H,1_B(\omega)v(y)\rangle
		=\langle h_B,v\rangle=0$
		for every $a$, $B$, and $v\in C(M_a)$.
		The density proved above gives $H=0$.
		Thus $\mathcal D^\perp=L^2(\mathcal F,\nu)$.
		
		\smallskip
		\noindent\emph{(iii).}
		Part~(i) gives
		$L^2(\mathcal F,\nu)\subseteq L^2(\W_s(S),\nu)$.
		If $\E(F\mid\W_s(S))=0$, then
		$F\in(\mathcal D^\perp)^\perp$.
		Since $(\mathcal D^\perp)^\perp$ is the
		$L^2(\nu)$-closure of $\mathcal D$,
		there exist $\psi_k\in\mathcal D$ with
		$\|F-\psi_k\|_{L^2(\nu)}\to0$.
	\end{proof}
	
	\begin{lemma}\label{lem:pair-orbit-measures}
		\label{lem:orbit-commutators}
		Let $\ell\ge3$, and let $(X,\mu,T_1,\ldots,T_\ell)$ be a system
		of invertible commuting measure-preserving transformations.
		Let $p_1,\ldots,p_\ell\in\Z[t]$ satisfy $p_i(0)=0$.
		Suppose that for each $i\in[\ell]$ there are
		$J_i,K_i\subseteq[\ell]$ with $i\in J_i\cap K_i$ and
		$J_i\cup K_i=[\ell]$, and that both $(p_j)_{j\in J_i}$ and
		$(p_j)_{j\in K_i}$ are linearly independent over $\Q$.
		For each $i$, let $\pi_i:(X,\mu,T_i)\to(Y_i,\nu_i,S_i)$ be a
		factor map onto a nilsystem bundle, where
		\[
		Y_i=\bigsqcup_{a\in I_i}\Omega_{i,a}\times M_{i,a},
		\qquad M_{i,a}=G_{i,a}/\Gamma_{i,a},
		\]
		and $S_i(\omega,y)=(\omega,\tau_{i,a}(\omega)y)$ on each piece.
		For $x\in X$, write
		$\pi_i(x)=(\omega_i,y_i)\in\Omega_{i,a_i}\times M_{i,a_i}$,
		and put $G_i=G_{i,a_i}$, $\Gamma_i=\Gamma_{i,a_i}$,
		$M_i=M_{i,a_i}$, $\tau_i=\tau_{i,a_i}(\omega_i)$, and
		$C_i=[G_i^\circ,G_i^\circ]$.
		
		Then there exists a measurable set $X_0\subseteq X$ with $\mu(X_0)=1$
		such that the following hold for every $x\in X_0$.
		\begin{enumerate}
			\item[(i)]
			The translation $y\mapsto\tau_i y$ is ergodic on $(M_i,m_{M_i})$,
			and $\pi_i(T_i^n x)=(\omega_i,\tau_i^n y_i)$ for every
			$i\in[\ell]$ and $n\in\Z$.
			\item[(ii)]
			The probability measures
			\[
			\lambda_{N,x}=\frac1N\sum_{n=1}^N
			\delta_{(\tau_1^{p_1(n)}y_1,\ldots,\tau_\ell^{p_\ell(n)}y_\ell)}
			\]
			converge weakly on $\prod_{i=1}^\ell M_i$ to a probability
			measure $\lambda_x$, where $\delta_z$ denotes the Dirac probability
			measure at $z$.
			\item[(iii)]
			For every $c_i\in C_i$, $i\in[\ell]$, and
			$\Phi\in C(\prod_{i=1}^\ell M_i)$,
			\[
			\int\Phi(c_1u_1,\ldots,c_\ell u_\ell)\,d\lambda_x(u)
			=\int\Phi(u_1,\ldots,u_\ell)\,d\lambda_x(u).
			\]
		\end{enumerate}
	\end{lemma}
	
	\begin{proof}
		By the bundle definition and the countability of the pairs $(i,n)$,
		there is a measurable set $X_*\subseteq X$ with $\mu(X_*)=1$
		such that, for every $x\in X_*$ and $i\in[\ell]$, the translation
		$y\mapsto\tau_i y$ is ergodic on $(M_i,m_{M_i})$ and
		$\pi_i(T_i^n x)=S_i^n\pi_i(x)$ for all $n\in\Z$.
		This proves (i) on $X_*$.
		
		For each $x\in X_*$, put $R_i=G_i^\circ y_i$.
		By \cite[Theorem~B and Section~1.9]{Leibman2005}, there are an
		integer $r\ge1$, a connected closed subgroup
		$H\le\prod_{i=1}^\ell G_i^\circ$, and points
		$z_b=(z_{b,1},\ldots,z_{b,\ell})\in\prod_{i=1}^\ell M_i$
		with compact orbits $Hz_b$, $1\le b\le r$, such that
		\[
		\lambda_{N,x}\longrightarrow
		\lambda_x=\frac1r\sum_{b=1}^r m_{Hz_b}
		\quad\text{weakly}.
		\]
		This proves (ii) on $X_*$.
		We shall prove (iii) on a measurable subset of $X_*$ of full measure.
		
		For every connected component $R$ of $M_{i,a}$, let $m_R$ be its normalized Haar
		measure.
		For fixed \(i,a\) and \(v\in C(M_{i,a})\), define two functions on \(X\)
		by
		\[
		\widehat v(x)=v(y),\qquad v_0(x)=\int_Rv\,dm_R
		\quad\text{if }\pi_i(x)=(\omega,y)\in\Omega_{i,a}\times R.
		\]
		Both functions are zero when \(\pi_i(x)\notin\Omega_{i,a}\times M_{i,a}\).
		The pair \(i,a\) is fixed whenever this notation is used.
		
		\smallskip
		\noindent\textbf{Claim 1.}
		\emph{For every such \(v\),}
		\[
		\E(\widehat v\mid\Krat(T_i))=v_0.
		\]
		
		Let \(\tau\) be an ergodic translation on \(M=G/\Gamma\), and let
		\(k\) be the number of its connected components. 
		Since $\tau$ is ergodic, no proper nonempty
		union of connected components can be
		$\tau$-invariant. Hence, for any connected component $R$,
		$M=\bigcup_{h=0}^{k-1}\tau^hR$. Thus \(\tau^kR=R\) and \(m_M(R)=1/k\) for every connected component \(R\).
		If \(D\subseteq R\) is \(\tau^k\)-invariant, then
		\[
		D'=\bigcup_{h=0}^{k-1}\tau^hD,\qquad
		m_M(\tau D'\mathbin{\triangle}D')=0,\qquad
		m_R(D)=m_M(D')\in\{0,1\}.
		\]
		Hence \(\tau^k|_R\) is ergodic. For \(R=G^\circ g\Gamma/\Gamma\),
		the open subgroup \(G_R=\langle G^\circ,\tau^k\rangle\) gives the
		identification
		\[
		G_R/(G_R\cap g\Gamma g^{-1})\longrightarrow R,\qquad
		h(G_R\cap g\Gamma g^{-1})\longmapsto hg\Gamma.
		\]
		Thus \(\tau^k|_R\) is a nilrotation on a connected space and is totally
		ergodic by \cite[Proposition~2.1(i)]{Frantzikinakis2008}.
		
		Let \(k\) be the number of components of \(M_{i,a}\), and fix \(m\ge k\).
		Since \(k\mid m!\), the preceding argument shows that
		\(\tau_{i,a}(\omega)^{m!}|_R\) is ergodic for almost every \(\omega\).
		For almost every $\omega\in\Omega_{i,a}$ and
		each connected component $R$, the mean ergodic
		theorem gives
		\[
		\begin{aligned}
			\int_R\left|
			\frac1N\sum_{n=1}^N
			v\bigl(\tau_{i,a}(\omega)^{m!n}y\bigr)
			-\int_Rv\,dm_R
			\right|^2dm_R(y)
			\longrightarrow0.
		\end{aligned}
		\]
		This integral is bounded by $4\|v\|_\infty^2$.
		Using $(\pi_i)_*\mu=\nu_i$,
		$m_{i,a}|_R=m_{i,a}(R)m_R$ and dominated
		convergence, we obtain

		\[
		\begin{aligned}
			&\left\|\frac1N\sum_{n=1}^N
			T_i^{m!n}\widehat v-v_0\right\|_{L^2(\mu)}^2\\
			&=\sum_R\int_{\pi_i^{-1}(\Omega_{i,a}\times R)}
			\left|\frac1N\sum_{n=1}^N T_i^{m!n}\widehat v(x)
			-\int_Rv\,dm_R\right|^2\,d\mu(x)\\
			&=\sum_R\int_{\Omega_{i,a}}\int_R
			\left|\frac1N\sum_{n=1}^N
			v\bigl(\tau_{i,a}(\omega)^{m!n}y\bigr)
			-\int_Rv\,dm_R\right|^2
			\,dm_{i,a}(y)\,d\beta_{i,a}(\omega)\\
			&=\sum_Rm_{i,a}(R)\int_{\Omega_{i,a}}\int_R
			\left|\frac1N\sum_{n=1}^N
			v\bigl(\tau_{i,a}(\omega)^{m!n}y\bigr)
			-\int_Rv\,dm_R\right|^2
			\,dm_R(y)\,d\beta_{i,a}(\omega)\\
			&\longrightarrow0.
		\end{aligned}
		\]
		Here \(R\) ranges over the \(k\) components.  Applying the mean ergodic theorem
		on \(X\) therefore gives
		\[
		\E(\widehat v\mid\I(T_i^{m!}))=v_0\qquad(m\ge k).
		\]
		Since \(\Krat(T_i)=\bigvee_{m\ge1}\I(T_i^{m!})\), convergence of
		increasing conditional expectations proves Claim~1.

		\smallskip
		\noindent\textbf{Claim 2.}
		\emph{There is a measurable set $X_0\subseteq X_*$ with $\mu(X_0)=1$
			such that, for every $x\in X_0$ and every nonempty
			$J\subseteq[\ell]$ for which $(p_i)_{i\in J}$ is linearly independent,}
		\[
		\prod_{i\in J}R_i\subseteq
		\operatorname{pr}_J(\operatorname{supp}\lambda_x),
		\qquad
		\operatorname{pr}_J(H)=\prod_{i\in J}G_i^\circ.
		\]
		Here $\operatorname{pr}_J$ denotes coordinate projection onto the
		indices in $J$.
		
		On each $M_{i,a}$ choose a countable base of nonempty open sets.
		For each basis set $U$, choose $v_U\in C(M_{i,a})$ with
		$0\le v_U\le1$, $v_U\not\equiv0$, and $\supp v_U\subseteq U$.
		Write $\mathcal T_{i,a}$ for this countable family.
		
		Fix a nonempty set $J$ for which $(p_i)_{i\in J}$ is linearly
		independent. For $i\in J$, fix $a_i\in I_i$ and
		$v_i\in\mathcal T_{i,a_i}$, and let $\widehat v_i,v_{i,0}$ be
		the functions from Claim~1 for the model $(i,a_i)$.
		Put
		\[
		q=\prod_{i\in J}|\operatorname{Comp}(M_{i,a_i})|,
		\qquad w=\prod_{i\in J}v_{i,0},
		\qquad B_N=\frac1N\sum_{n=1}^N
		\prod_{i\in J}T_i^{p_i(n)}\widehat v_i.
		\]
		Since $p_i(0)=0$, $q\mid p_i(qn)$.
		The component cycle calculation in Claim~1 and the identities on $X_*$
		give $v_{i,0}(T_i^{p_i(qn)}x)=v_{i,0}(x)$ for $x\in X_*$.
		It is clear that the polynomials remain linearly independent after replacing $t$ by $qt$ and then
		by Claim~1 and {\cite[Theorem~2.10]{FrantzikinakisKuca2025}},
		\[
		B_N^{(q)}:=\frac1N\sum_{n=1}^N
		\prod_{i\in J}T_i^{p_i(qn)}\widehat v_i
		\longrightarrow w\quad\text{in }L^2(\mu).
		\]
		Define the measurable set
		$E_J^{\mathbf a}(\mathbf v)
		=\{x\in X_*:w(x)>0,\ \lim_{N\to\infty}B_N(x)=0\}$,
		where $\mathbf a=(a_i)_{i\in J}$ and $\mathbf v=(v_i)_{i\in J}$.
		For every $x\in E_J^{\mathbf a}(\mathbf v)$, nonnegativity gives
		$0\le B_N^{(q)}(x)\le qB_{qN}(x)\to0$.
		Then by Fatou's lemma 
		\[
		\int_{E_J^{\mathbf a}(\mathbf v)}|w|^2\,d\mu
		\le\liminf_{N\to\infty}
		\|B_N^{(q)}-w\|_2^2=0.
		\]
		Since $w>0$ on $E_J^{\mathbf a}(\mathbf v)$, $\mu(E_J^{\mathbf a}(\mathbf v))=0$.
		Let
		\[
		\mathcal J=\{J\subseteq[\ell]:
		(p_i)_{i\in J}\text{ is linearly independent over }\Q,\ J\ne\varnothing\}.
		\]
		Hence $\mathcal J$ is finite, and then let
		
		\[
		X_0=X_*\setminus
		\bigcup_{\substack{J\in\mathcal J,\ a_i\in I_i\ (i\in J)\\
				v_i\in\mathcal T_{i,a_i}\ (i\in J)}}
		E_J^{\mathbf a}(\mathbf v).
		\]
		The union is countable, so $X_0$ is a measurable set with $\mu(X_0)=1$.
		
		Fix $x\in X_0$ and $J\in\mathcal J$.
		If the asserted support inclusion fails, there are nonempty open sets
		$U_i\subseteq R_i$, $i\in J$, such that
		\[
		\left(\prod_{i\in J}U_i\right)\cap
		\operatorname{pr}_J(\operatorname{supp}\lambda_x)=\varnothing.
		\]
		Choose $v_i\in\mathcal T_{i,a_i}$ with $\supp v_i\subseteq U_i$.
		The Haar measures $m_{R_i}$ give positive measure to every nonempty
		open set, so $w(x)=\prod_{i\in J}\int_{R_i}v_i\,dm_{R_i}>0$.
		On the other hand, the product $\prod_{i\in J}v_i(u_i)$ vanishes
		on $\operatorname{supp}\lambda_x$. Using $\pi_i(T_i^n x)=S_i^n\pi_i(x)$ and \(\lambda_{N,x}\to\lambda_x\),
		\[
		\begin{aligned}
			B_N(x)
			=\frac1N\sum_{n=1}^N
			\prod_{i\in J}v_i(\tau_i^{p_i(n)}y_i)
			&=\int\prod_{i\in J}v_i(u_i)\,d\lambda_{N,x}(u)\\
			&\longrightarrow
			\int\prod_{i\in J}v_i(u_i)\,d\lambda_x(u)=0.
		\end{aligned}
		\]
		Thus $x\in E_J^{\mathbf a}(\mathbf v)$, a contradiction.
		This proves the support inclusion.
		
		We now prove the group equality. Put
		$G_J=\prod_{i\in J}G_i^\circ$, $K=\operatorname{pr}_J(H)$,
		$R_J=\prod_{i\in J}R_i$, and $y_J=(y_i)_{i\in J}$.
		Since each Haar measure $m_{Hz_b}$ has full support
		on the compact orbit $Hz_b$,
		\[
		\operatorname{pr}_J(\operatorname{supp}\lambda_x)
		=\bigcup_{b=1}^r Kz_{b,J},
		\qquad z_{b,J}=(z_{b,i})_{i\in J}.
		\]
		
		Since $K$ is connected and $R_J$ is a product of
		connected components,
		\[
		Kz_{b,J}\cap R_J=
		\begin{cases}
			Kz_{b,J},&z_{b,J}\in R_J,\\
			\varnothing,&z_{b,J}\notin R_J.
		\end{cases}
		\]
		Then 
		$R_J
		=R_J\cap\bigcup_{b=1}^r Kz_{b,J}
		=\bigcup_{\substack{1\le b\le r\\z_{b,J}\in R_J}}
		Kz_{b,J}$.
		Note that $R_J=G_Jy_J$, so for any $b$ with $z_{b,J}\in R_J$
		we may choose $d_b\in G_J$ such that $d_by_J=z_{b,J}$.
		Let $D$ be the finite set of these elements. Then
		$R_J=\bigcup_{\xi\in D}K\xi y_J$.
		
		Write $y_J=(t_i\Gamma_i)_{i\in J}$. The stabilizer
		of $y_J$ in $G_J$ is
		\[
		\Lambda=\{g\in G_J:gy_J=y_J\}
		=\prod_{i\in J}(G_i^\circ\cap t_i\Gamma_i t_i^{-1}).
		\]
		It is discrete and countable.
		For $g\in G_J$, the preceding cover gives $gy_J=k\xi y_J$ for
		some $k\in K$ and $\xi\in D$. Thus $\xi^{-1}k^{-1}g\in\Lambda$,
		and $G_J=KD\Lambda$.
		
		Choose compact sets $F_m\subseteq H$ with $H=\bigcup_{m\ge1}F_m$.
		Then
		\[
		G_J=\bigcup_{\substack{m\ge1,\ \xi\in D\\\lambda\in\Lambda}}
		\operatorname{pr}_J(F_m)\xi\lambda.
		\]
		By Baire's theorem, one of these compact sets has nonempty interior. Then there exists a nonempty open set $U\subseteq K$.
		Since
		$1_{G_J}\in U^{-1}U
		=\bigcup_{u\in U}u^{-1}U\subseteq K$,
		the set $U^{-1}U$ is an open neighborhood of $1_{G_J}$.
		Hence
		\[
		K=\bigcup_{k\in K}kU^{-1}U,\qquad
		G_J\setminus K=\bigcup_{g\in G_J\setminus K}gK.
		\]
		Both sets are open. Since $G_J$ is connected, $K=G_J$. This completes Claim~2.
		
		\smallskip
		\noindent\textbf{Claim 3.}
		\emph{For every $x\in X_0$, the subgroup $H$ contains
			$\prod_{i=1}^\ell C_i$.}
		
		Fix $x\in X_0$, $i\in[\ell]$, and $a,b\in G_i^\circ$.
		Choose $J_i,K_i$ from the polynomial hypothesis.
		Claim~2 gives surjectivity of both $\operatorname{pr}_{J_i}(H)$
		and $\operatorname{pr}_{K_i}(H)$.
		We may therefore choose $u,v\in H$ such that
		\[
		\begin{aligned}
			u_i&=a,& u_j&=1_{G_j}\quad(j\in J_i\setminus\{i\}),\\
			v_i&=b,& v_j&=1_{G_j}\quad(j\in K_i\setminus\{i\}).
		\end{aligned}
		\]
		Since $J_i\cup K_i=[\ell]$, at least one of $u_j,v_j$ is the
		identity for each $j\ne i$. Consequently,
		\[
		[u,v]=(1_{G_1},\ldots,1_{G_{i-1}},[a,b],
		1_{G_{i+1}},\ldots,1_{G_\ell})\in H.
		\]
		These commutators generate $C_i$ in the $i$th coordinate.
		Taking products over $i$ gives $\prod_{i=1}^\ell C_i\subseteq H$.
		
		For $c_i\in C_i$ and $\Phi\in C(\prod_iM_i)$, by the invariance of Haar measure,
		\[
		\begin{aligned}
			\int\Phi(c_1u_1,\ldots,c_\ell u_\ell)\,d\lambda_x(u)
			&=\frac1r\sum_{b=1}^r
			\int_{Hz_b}\Phi(c_1u_1,\ldots,c_\ell u_\ell)\,dm_{Hz_b}(u)\\
			&=\frac1r\sum_{b=1}^r\int_{Hz_b}\Phi(u)\,dm_{Hz_b}(u)
			=\int\Phi(u)\,d\lambda_x(u).
		\end{aligned}
		\]
		This proves (iii). Since $X_0\subseteq X_*$, (i) and (ii) also
		hold on $X_0$. 
	\end{proof}
	
	\subsection{Proof of the phase-factor theorem}
	
	\begin{proof}[Integer assertion of Theorem~\ref{thm:phase-reduction}]
	Take $J_i,K_i$ from Lemma~\ref{lem:rank-cover}.
		For any distinct $i$ and $j$, at least one of $J_i$ or $K_i$ contains both indices. Consequently, $p_i$ and $p_j$ are linearly independent.
		Choose $s=s(D,\ell)\ge1$ from
		\cite[Theorem~2.8]{FrantzikinakisKuca2025}.
		We begin by assuming $E(f_1 \mid W_s(T_1)) = 0$ and proceed to show that this average tends to zero. By multilinearity, we may assume $|f_i| \leq 1$.
		
		Put $h_i=\E(f_i\mid\HK_s(T_i))$.
		By \cite[Theorem~1.2]{JamneshanMachado2026}, each $\HK_s(T_i)$
		is generated by an increasing sequence of $s$-step nilsystem bundle
		factors. For $\varepsilon>0$, choose this type of factor map
		$\pi_i:(X,\mu,T_i)\to(Y_i,\nu_i,S_i)$  such that
		\[
		\begin{gathered}
			\pi_i^{-1}\mathcal Y_i\subseteq\HK_s(T_i),\qquad
			F_i=\E(f_i\mid\pi_i^{-1}\mathcal Y_i),\\
			\sum_{i=1}^\ell\|h_i-F_i\|_2<\varepsilon,
		\end{gathered}
		\]
		where $\mathcal Y_i$ is the sigma-algebra of $Y_i$.
		We identify functions on $Y_i$ with their pullbacks to $X$.
		The choice of $s$ gives
		$\|A_N^{\mathbf p}(f_1,\ldots,f_\ell)
		-A_N^{\mathbf p}(h_1,\ldots,h_\ell)\|_2\to0$.
		Since $|h_i|,|F_i|\le1$, expand the differences coordinate by coordinate, and using the property of measure preservation, obtain
		\[
		\|A_N^{\mathbf p}(h_1,\ldots,h_\ell)
		-A_N^{\mathbf p}(F_1,\ldots,F_\ell)\|_2
		\le\sum_{i=1}^\ell\|h_i-F_i\|_2.
		\]
		Consequently,
		\[
		\limsup_{N\to\infty}
		\|A_N^{\mathbf p}(f_1,\ldots,f_\ell)
		-A_N^{\mathbf p}(F_1,\ldots,F_\ell)\|_2<\varepsilon.
		\]
		It remains to prove that
		$\|A_N^{\mathbf p}(F_1,\ldots,F_\ell)\|_2\to0$.
		
		Let $\mathcal V_i,\mathcal D_i$ be the spaces in
		Lemma~\ref{lem:commutator-density} for $Y_i$.
		For a unit-modulus function $\phi$ on $Y_1$ and every $q\ge1$,
		\[
		\Delta_{T_1^q}^{s+1}(\phi\circ\pi_1)
		=(\Delta_{S_1^q}^{s+1}\phi)\circ\pi_1.
		\]
		Hence $\pi_1^{-1}\W_s(S_1)\subseteq\W_s(T_1)$.
		For $h\in L^2(\W_s(S_1),\nu_1)$, it follows that
		\[
		\int_{Y_1}F_1\overline h\,d\nu_1
		=\int_X f_1\,\overline{h\circ\pi_1}\,d\mu=0.
		\]
		By Lemma~\ref{lem:commutator-density}(iii), there are
		$\psi_k\in\mathcal D_1$ with $\|F_1-\psi_k\|_{L^2(\nu_1)}\to0$.
		
		We first prove $L^2$ convergence to zero for
		$V_1\in\mathcal D_1$ and $V_i\in\mathcal V_i$, $2\le i\le\ell$.
		By multilinearity and the definitions of these spaces, it suffices
		to take, on fixed pieces $\Omega_{i,a_i}\times M_{i,a_i}$,
		\[
		V_1(\omega,y)=b_1(\omega)(v_1(cy)-v_1(y)),
		\qquad
		V_i(\omega,y)=b_i(\omega)v_i(y)\quad(2\le i\le\ell),
		\]
		where $b_i\in L^\infty(\beta_{i,a_i})$,
		$v_i\in C(M_{i,a_i})$, and
		$c\in[G_{1,a_1}^\circ,G_{1,a_1}^\circ]$. For all $i$, set $V_i=0$ outside
		$\Omega_{i,a_i}\times M_{i,a_i}$.
		
		Let $X_0$ be the set in Lemma~\ref{lem:pair-orbit-measures}.
		If $\pi_i(x)\notin\Omega_{i,a_i}\times M_{i,a_i}$ for some $i$
		and $x\in X_0$, then
		$A_N^{\mathbf p}(V_1,\ldots,V_\ell)(x)=0$ for every $N$.
    	We may therefore assume that
    	$\pi_i(x)=(\omega_i,y_i)\in\Omega_{i,a_i}\times M_{i,a_i}$
	    for every $i$.
		By Lemma~\ref{lem:pair-orbit-measures} (ii) and (iii) we have
		\[
		\begin{aligned}
			\lim_{N\to\infty}A_N^{\mathbf p}(V_1,\ldots,V_\ell)(x)
			&=\left(\prod_{i=1}^\ell b_i(\omega_i)\right)
			\int (v_1(cu_1)-v_1(u_1))
			\prod_{i=2}^\ell v_i(u_i)\,d\lambda_x(u)\\
			&=0.
		\end{aligned}
		\]
		The last equality uses invariance under $(c,1_{G_2},\ldots,1_{G_\ell})$. Using
		multilinearity and dominated convergence one can obtain the required $L^2$ limit
		for all $V_1\in\mathcal D_1$ and $V_i\in\mathcal V_i$, $i\ge2$.
		
		For any $g_1,\ldots,g_\ell\in L^\infty(\mu)$ and every $j\in[\ell]$,
		measure preservation gives
		\[
		\sup_N\|A_N^{\mathbf p}(g_1,\ldots,g_\ell)\|_2
		\le\|g_j\|_2\prod_{i\ne j}\|g_i\|_\infty.
		\]
		Fix $\psi\in\mathcal D_1$ and apply Lemma~\ref{lem:l2extension}
		to the remaining $\ell-1$ coordinates. Since
		the spaces $\mathcal V_i$ are dense in $L^2(\nu_i)$, 
		$A_N^{\mathbf p}(\psi,F_2,\ldots,F_\ell)\to0$ in $L^2$.
		Taking $\psi=\psi_k$ and using $|F_i|\le1$ yields
		\[
		\limsup_{N\to\infty}
		\|A_N^{\mathbf p}(F_1,\ldots,F_\ell)\|_2
		\le\|F_1-\psi_k\|_2\longrightarrow0.
		\]
		By the preceding estimate we have
		$\limsup_{N\to\infty}\|A_N^{\mathbf p}(f_1,\ldots,f_\ell)\|_2\le\varepsilon$.
		Since $\varepsilon>0$ is arbitrary, the limit is zero.
		The same argument applies to any coordinate $j$ with
		$\E(f_j\mid\W_s(T_j))=0$.
		
		Finally let $g_i=\E(f_i\mid\W_s(T_i))$.
		Then $\E(f_i-g_i\mid\W_s(T_i))=0$, and by multilinearity we have
		\[
		\begin{aligned}
			A_N^{\mathbf p}(f_1,\ldots,f_\ell)
			-A_N^{\mathbf p}(g_1,\ldots,g_\ell)
			=\sum_{j=1}^\ell
			A_N^{\mathbf p}(g_1,\ldots,g_{j-1},f_j-g_j,f_{j+1},\ldots,f_\ell).
		\end{aligned}
		\]
		Each term tends to zero in $L^2(\mu)$, as proved above.
		This proves the integer characteristic-factor assertion.
	\end{proof}
	
\section{Degree bounds for polynomial phases}\label{sec:interpolation}

In this section, we show the
degree bounds $4,2,2$ for the family $t,t^2,t^2+t$ and then prove the integer upper bounds in Theorem~\ref{thm:main}.

	Recall the coefficientwise reduction $P\mapsto[P]$ modulo $\Q$
	from Section~\ref{sec:preliminaries}. We need the Weyl's theorem.
	
	\begin{lemma}[Weyl's theorem]\label{lem:weyl}
		Let $m\ge0$ be an integer and
		$P(t)=a_0+a_1t+\cdots+a_mt^m\in\R[t]$.
		If
		\[
		\limsup_{N\to\infty}
		\left|\frac1N\sum_{n=1}^N e(P(n))\right|>0,
		\]
		then $a_k\in\Q$ for every $1\le k\le m$.
		Equivalently, the polynomial $[P]\in(\R/\Q)[t]$ obtained
		by reducing the coefficients of $P$ modulo $\Q$ is constant.
	\end{lemma}
	
	\begin{proof}
		If $a_k$ is irrational for some $k\ge1$, Weyl's
		equidistribution theorem \cite{Weyl1916} gives
		\[
		\frac1N\sum_{n=1}^N e(P(n))\longrightarrow0,
		\]
		contrary to the assumption.
	\end{proof}
	
	We apply this criterion to compute conditional expectations on phase
	factors. 	
	\begin{lemma}\label{lem:phaseprojection}
		Let $(X,\mu,S)$ be a measure-preserving system, let $d,L\ge1$
		be integers, and let
		$\phi\in L^\infty(\mu)$ satisfy
		$|\phi|=1$ and $\Delta_S^L\phi=1$.
		Put
		\[
		E=\bigcap_{k=d+1}^{L-1}
		\{x:\Delta_S^k\phi(x)\in\roots\},
		\]
		with $E=X$ when $d\ge L-1$. Then $E$ is $S$-invariant and
		$\E(\phi\mid\W_d(S))=1_E\phi$.
	\end{lemma}
	
	\begin{proof}
		If $d\ge L-1$, then $E=X$ and $\phi\in\mathcal P_d(S)$,
		so the conclusion holds. Suppose that $d<L-1$.
		Put $\psi_k=\Delta_S^k\phi$ and, for every $r\ge1$, define
		$E_r=X_{S,\phi}\cap\bigcap_{k=d+1}^{L-1}
		\{x:\psi_k(x)^{r!}=1\}$,
		where $X_{S,\phi}$ is the invariant conull set from
		Lemma~\ref{lem:polynomial-phases}.
		For every $x\in X_{S,\phi}$, we have
		$\psi_k(Sx)=\psi_k(x)\psi_{k+1}(x)$ and $\psi_L(x)=1$.
		Descending induction on $k$ gives $S^{-1}E_r=E_r$.
		Since $r!\mid(r+1)!$, we have $E_r\subseteq E_{r+1}$.
		Moreover, every root of unity has finite order, and there
		are only finitely many $k$. Hence
		$E\cap X_{S,\phi}=\bigcup_{r\ge1}E_r$.
		Since $\mu(X_{S,\phi})=1$, $E$ is $S$-invariant.
		
		We first prove that $1_E\phi$ is $\W_d(S)$-measurable.
		Fix $r\ge1$ and put $M=r!L!$. For $n\in\Z$ and $1\le k<L$,
		\[
		\binom{Mn}{k}
		=r!n\,\frac{L!}{k!}\prod_{b=1}^{k-1}(Mn-b)\in r!\Z.
		\]
		Hence, for $x\in E_r$, by Lemma~\ref{lem:polynomial-phases}(i),
		$\phi(S^{Mn}x)=\prod_{k=0}^{d}\psi_k(x)^{\binom{Mn}{k}}$.
		The exponents have degree at most $d$ in $n$, so
		$\Delta_{S^M}^{d+1}\phi=1$ on $E_r$.
		Since $E_r$ is invariant,
		$1_{E_r}\phi+1_{E_r^c}\in\mathcal P_d(S)$.
		Also $E_r\in\I(S)\subseteq\W_d(S)$, so
		$1_{E_r}\phi\in L^\infty(\W_d(S))$.
		As $1_{E_r}\phi\to1_E\phi$ almost everywhere, we obtain
		$1_E\phi\in L^\infty(\W_d(S))$.
		
		It remains to show that $1_{E^c}\phi$ is orthogonal to
		$L^2(\W_d(S))$. By Lemma~\ref{lem:phasefactors}(iii),
		it suffices to take $\chi\in\mathcal P_d(S)$.
		Choose $q\ge1$ with $\Delta_{S^q}^{d+1}\chi=1$.
		For almost every $x\in E^c$, Lemma~\ref{lem:polynomial-phases}(ii)
		gives $\deg[P_{S,\phi,x}]>d$, while
		$\deg P_{S^q,\chi,x}\le d$. Therefore
		\[
		\deg\bigl[P_{S,\phi,x}(qt)-P_{S^q,\chi,x}(t)\bigr]
		=\deg[P_{S,\phi,x}]>d.
		\]
		By Lemma~\ref{lem:polynomial-phases} (i) and Weyl's theorem \cite{Weyl1916},
		\[
		\begin{aligned}
			\frac1N\sum_{n=1}^N
			\phi(S^{qn}x)\overline{\chi(S^{qn}x)}
			=\frac1N\sum_{n=1}^N
			e\bigl(P_{S,\phi,x}(qn)-P_{S^q,\chi,x}(n)\bigr)
			\longrightarrow0.
		\end{aligned}
		\]
		Since $E$ is invariant and $|\phi\overline\chi|=1$,
		by dominated convergence, we have
		\[
		\begin{aligned}
			\int_X1_{E^c}\phi\overline\chi\,d\mu
			=\lim_{N\to\infty}\int_X1_{E^c}(x)
			\frac1N\sum_{n=1}^N
			\phi(S^{qn}x)\overline{\chi(S^{qn}x)}\,d\mu(x)
			=0.
		\end{aligned}
		\]
		Thus $1_{E^c}\phi\perp L^2(\W_d(S))$. Together with the
		measurability of $1_E\phi$, this proves the formula.
	\end{proof}
	
	\begin{lemma}\label{lem:polynomial}
		Let $V$ be a vector space over $\Q$, let $\delta\in\Q\setminus\{0\}$,
		and let $P_0,P_\delta\in V[u]$ and $Q\in V[u,w]$ be formal
		polynomials with coefficients in $V$, with $\deg0=-\infty$.
		Suppose that, for $\beta\in\{0,\delta\}$, the polynomial
		\[
		P_\beta(u+t)+Q(u,\beta+t^2+t)
		\]
		is unchanged when $t$ is replaced by $-t$. Then
		\[
		\deg P_0\le4,\qquad \deg P_\delta\le4,\qquad \deg_wQ\le2.
		\]
	\end{lemma}
	
	\begin{proof}
		For $\beta\in\{0,\delta\}$, put
		\[
		\begin{aligned}
			H_\beta(u,t)&=P_\beta(u+t)+Q(u,\beta+t^2+t),\\
			R_\beta(u,w)&=(4(w-\beta)+1)Q_{ww}(u,w)
			+2Q_w(u,w)-Q_{uu}(u,w).
		\end{aligned}
		\]
		The chain rule gives
		\[
		(\partial_t^2-\partial_u^2)H_\beta(u,t)
		=R_\beta(u,\beta+t^2+t).
		\]
		The left-hand side is even in $t$. If $R_\beta$ had positive
		$w$-degree $m$ with leading coefficient $a_m(u)\ne0$, the coefficient
		of $t^{2m-1}$ on the right would be $ma_m(u)\ne0$; all lower
		$w$-degree terms have $t$-degree at most $2m-2$.
		Thus $R_0$ and $R_\delta$ are independent of $w$. Since
		\[
		R_\delta(u,w)-R_0(u,w)=-4\delta Q_{ww}(u,w),
		\]
		$Q_{ww}$ is independent of $w$, and $\deg_wQ\le2$.
		Consequently, $Q(u,\beta+t^2+t)$ has no $t^5$ term.
		Taylor's formula shows that the coefficient of $t^5$ in
		$H_\beta(u,t)$ is $P_\beta^{(5)}(u)/5!$.
		Evenness gives $P_\beta^{(5)}=0$, so $\deg P_\beta\le4$.
	\end{proof}
	
	\begin{proposition}\label{prop:phases}
		Let $T_1,T_2,T_3$ be invertible commuting measure-preserving
		transformations of a standard probability space $(X,\mathcal X,\mu)$.
		Suppose that $\phi_1,\phi_2,\phi_3\in L^\infty(\mu)$ and that
		there are integers $q,L\ge1$ such that, for each $i=1,2,3$,
		\[
		|\phi_i|=1,\qquad \Delta_{T_i^q}^L\phi_i=1.
		\]
		Let $A_N$ be the average in \eqref{eq:average}.
		Put $(d_1,d_2,d_3)=(4,2,2)$ and, for each $i=1,2,3$, define
		\[
		h_i=\phi_i-\E(\phi_i\mid\W_{d_i}(T_i)).
		\]
		Then
		\[
		\begin{aligned}
			\|A_N(h_1,\phi_2,\phi_3)\|_2&\longrightarrow0,\\
			\|A_N(\phi_1,h_2,\phi_3)\|_2&\longrightarrow0,\\
			\|A_N(\phi_1,\phi_2,h_3)\|_2&\longrightarrow0.
		\end{aligned}
		\]
	\end{proposition}
	
	\begin{proof}
		Fix $j\in\{1,2,3\}$.
		Apply Lemma~\ref{lem:polynomial-phases} to each $(T_i^q,\phi_i)$. Set
		\[
		\begin{aligned}
			X_0=\bigcap_{g\in\langle T_1,T_2,T_3\rangle}g^{-1}
			\Bigl(\bigcap_{i=1}^3X_{T_i^q,\phi_i}
			\cap\bigcap_{i,k=1}^3\{y\in X:T_iT_ky=T_kT_iy\}\Bigr).
		\end{aligned}
		\]
		Then $\mu(X_0)=1$ and $gX_0=X_0$ for every
		$g\in\langle T_1,T_2,T_3\rangle$.
		By Lemma~\ref{lem:phasefactors},
		$\W_{d_j}(T_j^q)=\W_{d_j}(T_j)$. Put $V=\R/\Q$, regarded as a $\Q$-vector space. For every
		$y\in X_0$, write $P_{i,y}=P_{T_i^q,\phi_i,y}$ for the polynomials
		from Lemma~\ref{lem:polynomial-phases}.
		Lemma~\ref{lem:polynomial-phases}(ii) and
		Lemma~\ref{lem:phaseprojection} therefore give
		\[
		\E(\phi_j\mid\W_{d_j}(T_j))(y)
		=1_{\{\deg[P_{j,y}]\le d_j\}}\phi_j(y)
		\quad\text{for almost every }y.
		\]
		Let 
		$h_j(y)=1_{\{\deg[P_{j,y}]>d_j\}}\phi_j(y)$ for every $y\in X_0$
		and set $h_j=0$ outside $X_0$.
		
		Write $(p_1(t),p_2(t),p_3(t))=(t,t^2,t^2+t)$.
		Fix an integer $0\le a<q$, put $b_i=p_i(-a)$, and define
		\[
		Q_{i,y}(t)=\left[P_{i,T_i^{b_i}y}
		\left(\frac{t-b_i}{q}\right)\right]\in V[t].
		\]
		For $h\in q\Z$, Lemma~\ref{lem:polynomial-phases}(ii) gives
		$\deg Q_{i,y}<L$, and $Q_{i,T_i^hy}(t)=Q_{i,y}(t+h)$.
		Since $b_j=p_j(-a)$ and $p_j\in\Z[t]$, we have
		\[
		m:=\frac{p_j(qn-a)-b_j}{q}\in\Z.
		\]
		Thus
		$T_j^{p_j(qn-a)}y=T_j^{qm}T_j^{b_j}y$.
		By Lemma~\ref{lem:polynomial-phases}(ii),
		$
		[P_{j,T_j^{p_j(qn-a)}y}](t)
		=[P_{j,T_j^{b_j}y}](t+m).
		$
		By the definition of $Q_{j,y}$, 
		\[
		\deg[P_{j,T_j^{p_j(qn-a)}y}]
		=\deg[P_{j,T_j^{b_j}y}]
		=\deg Q_{j,y}.
		\]
		Consequently,
		\[
		h_j(T_j^{p_j(qn-a)}y)
		=1_{\{\deg Q_{j,y}>d_j\}}\phi_j(T_j^{p_j(qn-a)}y).
		\]
		Put
		\[
		B_{M,a}(y)=\frac1M\sum_{n=1}^M
		h_j(T_j^{p_j(qn-a)}y)
		\prod_{i\ne j}\phi_i(T_i^{p_i(qn-a)}y).
		\]
		
		\smallskip
		\noindent\textbf{Claim.}
		\emph{The measurable set}
		\[
		E=\{y\in X_0:\limsup_{M\to\infty}|B_{M,a}(y)|>0\}
		\]
		\emph{has measure zero.}
		
		Put
		$
		F_y(n)=\sum_{i=1}^3P_{i,T_i^{b_i}y}
		\left(\frac{p_i(qn-a)-b_i}{q}\right).
		$
		For $y\in E$, we have $\deg Q_{j,y}>d_j$. Hence
		\[
		\begin{aligned}
			B_{M,a}(y)
			=\frac1M\sum_{n=1}^M
			\prod_{i=1}^3\phi_i(T_i^{p_i(qn-a)}y)
			=\frac1M\sum_{n=1}^M e(F_y(n)),
		\end{aligned}
		\]
		where the last equality follows from
		Lemma~\ref{lem:polynomial-phases}(i).
		By the definition of $E$,
		\[
		\limsup_{M\to\infty}
		\left|\frac1M\sum_{n=1}^M e(F_y(n))\right|>0.
		\]
		Lemma~\ref{lem:weyl} therefore gives that $[F_y]$ is constant.	
		Substituting $n=(t+a)/q$, we obtain
		\[
		\begin{aligned}
			\left[F_y\left(\frac{t+a}{q}\right)\right]
			&=\sum_{i=1}^3
			\left[P_{i,T_i^{b_i}y}
			\left(\frac{p_i(t)-b_i}{q}\right)\right]
			=\sum_{i=1}^3Q_{i,y}(p_i(t))\\
			&=Q_{1,y}(t)+Q_{2,y}(t^2)+Q_{3,y}(t^2+t).
		\end{aligned}
		\]
		The left-hand side is constant in $t$, since $[F_y]$
		is constant. Taking $t=0$ gives
		\[
		Q_{1,y}(t)+Q_{2,y}(t^2)+Q_{3,y}(t^2+t)
		=\sum_{i=1}^3Q_{i,y}(0).
		\]
		Then
		$Q_{1,y}(t)+Q_{3,y}(t^2+t)$ is even in $t$.
		
		Suppose that $\mu(E)>0$. Recurrence for $T_3^q$ gives $\delta\in q\N$
		with $\mu(E\cap T_3^{-\delta}E)>0$.
		Recurrence for $T_1^q$ on $E\cap T_3^{-\delta}E$ gives a point
		$x\in E\cap T_3^{-\delta}E$ and $L$ distinct return times
		\[
		0=\alpha_0<\alpha_1<\cdots<\alpha_{L-1},\qquad
		\alpha_r\in q\Z_{\ge0},
		\]
		such that, $x_r:=T_1^{\alpha_r}x\in E\cap T_3^{-\delta}E$
		for every $0\le r<L$.
		For $0\le r<L$, define
		\[
		\ell_r(u)=\prod_{\substack{0\le k<L\\k\ne r}}
		\frac{u-\alpha_k}{\alpha_r-\alpha_k},\qquad
		Q(u,w)=\sum_{r=0}^{L-1}\ell_r(u)Q_{3,x_r}(w).
		\]
		For $0\le r,s<L$,
		\[
		\ell_r(\alpha_s)=
		\begin{cases}
			1,&r=s,\\
			0,&r\ne s,
		\end{cases}
		\qquad
		Q(\alpha_s,w)=Q_{3,x_s}(w).
		\]
		In particular, by $\alpha_0=0$ and $x_0=x$, one has $Q(0,w)=Q_{3,x}(w)$.
		
		For $\beta\in\{0,\delta\}$, the identity
		$Q_{i,T_i^hy}(t)=Q_{i,y}(t+h)$ $(h\in q\Z)$
		and $T_3^\beta x_r=T_1^{\alpha_r}T_3^\beta x$, by commutativity,
		gives
		\[
		\begin{aligned}
			Q_{1,T_3^\beta x_r}(t)=Q_{1,T_3^\beta x}(\alpha_r+t),\quad
			Q_{3,T_3^\beta x_r}(w)=Q_{3,x_r}(w+\beta).
		\end{aligned}
		\]
		Note that the Lagrange interpolation formula shows that a polynomial of degree less than $L$ with coefficients in a
		$\Q$-vector space is determined by its values at $L$ distinct
		rational points.
		Since $Q_{1,T_3^\beta x}(u+t)$ has degree less than $L$ in $u$,
		Lagrange interpolation gives
		\[
		\begin{aligned}
			&Q_{1,T_3^\beta x}(u+t)+Q(u,\beta+t^2+t)\\
			&\quad=\sum_{r=0}^{L-1}\ell_r(u)
			\bigl(Q_{1,T_3^\beta x_r}(t)
			+Q_{3,T_3^\beta x_r}(t^2+t)\bigr).
		\end{aligned}
		\]
		Each bracketed polynomial is even in $t$, since $T_3^\beta x_r\in E$.
		The coefficients $\ell_r(u)$ do not depend on $t$, so the left-hand
		side is even for both $\beta=0$ and $\beta=\delta$.
		Applying Lemma~\ref{lem:polynomial} with
		$P_0=Q_{1,x}$ and $P_\delta=Q_{1,T_3^\delta x}$ gives
		\[
		\deg Q_{1,x}\le4,\qquad \deg_wQ\le2.
		\]
		Since $Q(0,w)=Q_{3,x}(w)$, we have $\deg Q_{3,x}\le2$.
		Also, the identity at $x$ gives
		\[
		Q_{2,x}(t^2)=\sum_{i=1}^3Q_{i,x}(0)
		-Q_{1,x}(t)-Q_{3,x}(t^2+t).
		\]
		Its right-hand side has degree at most four, so $\deg Q_{2,x}\le2$.
		Thus $\deg Q_{j,x}\le d_j$, contrary to $x\in E$.
		This proves the claim.
		
		For each $0\le a<q$, the claim gives
$B_{M,a}(y)\longrightarrow0$
		for almost every $y$.
		Since there are only finitely many values of $a$,
		there is a conull set on which these limits hold for all $a$.
		
		Let
		\[
		J_N(y)=\frac1N\sum_{n=1}^N
		h_j(T_j^{p_j(n)}y)
		\prod_{i\ne j}\phi_i(T_i^{p_i(n)}y).
		\]
		For every $M\ge1$, the integers $qn-a$, where
		$1\le n\le M$ and $0\le a<q$, are 
		$1,\ldots,qM$. Hence
		\[
		J_{qM}(y)=\frac1q\sum_{a=0}^{q-1}B_{M,a}(y)
		\longrightarrow0
		\quad\text{for almost every }y.
		\]
		
		For $N\ge q$, put $M=\lfloor N/q\rfloor$, so that
		$0\le N-qM<q$. Then
		\[
		\begin{aligned}
			J_N(y)
			=\frac{qM}{N}J_{qM}(y)
			+\frac1N\sum_{n=qM+1}^N
			h_j(T_j^{p_j(n)}y)
			\prod_{i\ne j}\phi_i(T_i^{p_i(n)}y).
		\end{aligned}
		\]
		Since $|h_j|\le1$ and $|\phi_i|=1$,
		\[
		|J_N(y)|
		\le\frac{qM}{N}|J_{qM}(y)|+\frac{N-qM}{N}
		\le |J_{qM}(y)|+\frac qN
		\]
		for almost every $y$.
		Thus $J_N(y)\to0$ almost everywhere.
		Finally, $|J_N|\le1$, so dominated convergence gives
		\[
		\lim_{N\to\infty}\|J_N\|_2^2
		=\lim_{N\to\infty}\int_X|J_N|^2\,d\mu=0.
		\]
		This completes the proof.
	\end{proof}
	
	\subsection{Proof of the sharp-order theorem}\label{sec:mainproof}
		Now we shall prove the integer upper bounds in Theorem~\ref{thm:main}. The prime assertions
	are proved in Section~\ref{sec:primes}.

	\begin{proof}[Integer upper bounds in Theorem~\ref{thm:main}
		and Corollary~\ref{cor:orders}]
		Put $(d_1,d_2,d_3)=(4,2,2)$ and
		$(p_1(t),p_2(t),p_3(t))=(t,t^2,t^2+t)$.
		Fix $j\in\{1,2,3\}$ and suppose that
		$\E(f_j\mid\W_{d_j}(T_j))=0$.
		By multilinearity, we may assume $|f_i|\le1$ for $i=1,2,3$.
		The polynomials $t,t^2,t^2+t$ are pairwise linearly independent.
		Apply Corollary~\ref{cor:half-independent} with $\ell=3$ and $D=2$,
		increase $s$ to at least four if necessary, and put
		$F_i=\E(f_i\mid\W_s(T_i))$.
		By Corollary~\ref{cor:half-independent} and the tower property of conditional expectation, we have $\E(F_j\mid\W_{d_j}(T_j))
		=\E(f_j\mid\W_{d_j}(T_j))=0$ and
		\[
		\begin{aligned}
			\|A_N(f_1,f_2,f_3)-A_N(F_1,F_2,F_3)\|_2
			\longrightarrow0.
		\end{aligned}
		\]
		Next we will prove that $\|A_N(F_1,F_2,F_3)\|_2\to0$.
		
		Recall that $\mathcal P_d(T):=\{\phi\in L^\infty(\mu):|\phi|=1,
		\ \Delta_{T^q}^{d+1}\phi=1\text{ for some }q\ge1\}$.
		Fix arbitrary $U_i\in\mathcal A_s(T_i):=\operatorname{span}_{\C}\mathcal P_s(T_i)$
		and define $V_j=U_j-\E(U_j\mid\W_{d_j}(T_j))$,
		$V_i=U_i$ for $i\ne j$.
		Write each $U_i$ as a finite sum
		$U_i=\sum_{r=1}^{m_i}c_{i,r}\phi_{i,r}$, where
		$c_{i,r}\in\C$ and $\phi_{i,r}\in\mathcal P_s(T_i)$.
		For each $i=1,2,3$, choose $q_i\ge1$ such that
		$\Delta_{T_i^{q_i}}^{s+1}\phi_{i,r}=1$
		for every $1\le r\le m_i$.
		Let $q$ be the least common multiple of $q_1,q_2,q_3$.
		Since $q_i\mid q$, Lemma~\ref{lem:phasefactors} gives
		$\Delta_{T_i^q}^{s+1}\phi_{i,r}=1$
		for all $i,r$.
		Proposition~\ref{prop:phases}, applied with $L=s+1$
		to each triple of phases, and multilinearity therefore give
		\[
		\|A_N(V_1,V_2,V_3)\|_2\longrightarrow0
		\]
		for every $U_i\in\mathcal A_s(T_i)$.
		
		Since $\E(F_j\mid\W_{d_j}(T_j))=0$, by
		the $L^2$ contraction property of conditional expectation, we have
		\[
		\begin{aligned}
			\|F_j-V_j\|_2
			=\|(F_j-U_j)
			-\E(F_j-U_j\mid\W_{d_j}(T_j))\|_2
			\le\|F_j-U_j\|_2.
		\end{aligned}
		\]
		For $i\ne j$, we have $F_i-V_i=F_i-U_i$. Hence, for $i=1,2,3$,
$\|F_i-V_i\|_2\le\|F_i-U_i\|_2$.

		Now fix $\varepsilon>0$. By Lemma~\ref{lem:phasefactors},
		first choose $U_3$ and define $V_3$, then choose $U_2$
		and define $V_2$, and finally choose $U_1$, so that
		\[
		\begin{aligned}
			\|F_3-U_3\|_2
			&<\frac{\varepsilon}{3},\\
			\|F_2-U_2\|_2
			&<\frac{\varepsilon}{3(1+\|V_3\|_\infty)},\\
			\|F_1-U_1\|_2
			&<\frac{\varepsilon}
			{3(1+\|V_2\|_\infty)(1+\|V_3\|_\infty)}.
		\end{aligned}
		\]
		For any $a_1,a_2,a_3\in L^\infty(\mu)$, by measure preservation, we have
		\[
		\begin{aligned}
			\|A_N(a_1,a_2,a_3)\|_2
			&\le\frac1N\sum_{n=1}^N
			\|T_i^{p_i(n)}a_i\|_2
			\prod_{k\ne i}\|a_k\|_\infty\\
			&=\|a_i\|_2\prod_{k\ne i}\|a_k\|_\infty,
			\qquad \text{for~}i=1,2,3.
		\end{aligned}
		\]
		Since $|F_i|\le1$, by multilinearity and the preceding estimate 
		\[
		\begin{aligned}
			\|A_N(F_1,F_2,F_3)-A_N(V_1,V_2,V_3)\|_2
			&\le
			\|F_3-V_3\|_2
			+\|F_2-V_2\|_2\|V_3\|_\infty\\
			&\qquad+
			\|F_1-V_1\|_2
			\|V_2\|_\infty\|V_3\|_\infty\\
			&<\varepsilon.
		\end{aligned}
		\]
		Since $\varepsilon>0$ is arbitrary and $\|A_N(V_1,V_2,V_3)\|_2\longrightarrow0$, then
		$\|A_N(F_1,F_2,F_3)\|_2\longrightarrow0$.
		Hence
		\[
		\|A_N(f_1,f_2,f_3)\|_2\longrightarrow0.
		\]
		
		By Lemma~\ref{lem:phasefactors}, we have
		$\W_{d_j}(T_j)\subseteq\HK_{d_j}(T_j)$.
		Thus
		$\E(f_j\mid\HK_{d_j}(T_j))=0$ implies $\E(f_j\mid\W_{d_j}(T_j))=0$.
		Together with
		$\|f_j\|_{U^{d_j+1}(T_j)}=0$ if and only if $
		\E(f_j\mid\HK_{d_j}(T_j))=0$,
		this proves the seminorm assertions.
		
		Finally, for every $i=1,2,3$, let
		$g_i=\E(f_i\mid\W_{d_i}(T_i))$ or $g_i=\E(f_i\mid\HK_{d_i}(T_i))$.
		In either case, by the  property of conditional expectation,
		$\E(f_i-g_i\mid\W_{d_i}(T_i))=0$.
		By multilinearity,
		\[
		\begin{aligned}
			A_N(f_1,f_2,f_3)-A_N(g_1,g_2,g_3)
			&=A_N(f_1-g_1,f_2,f_3)
			+A_N(g_1,f_2-g_2,f_3)\\
			&\qquad+A_N(g_1,g_2,f_3-g_3).
		\end{aligned}
		\]
		By the previously proved result that the limit is zero, each term converges to zero in $L^2(\mu)$. Hence
		\[
		\|A_N(f_1,f_2,f_3)-A_N(g_1,g_2,g_3)\|_2
		\longrightarrow0.
		\]
	\end{proof}

	\section{A five-dimensional sharpness example}\label{sec:sharpness}
	
	\begin{proposition}\label{prop:sharpness}
		There exist commuting Haar measure-preserving affine automorphisms \(T,S\)
		of \(\T^5\), equipped with Haar probability \(m_{\T^5}\), and
 \(f_1,f_2,f_3\in L^\infty(m_{\T^5})\) such that
		\[
		\norm{f_1}_{U^4(T)}=0,\qquad
		\norm{f_2}_{U^2(S)}=\norm{f_3}_{U^2(S)}=0,
		\]
		whereas for all $N\ge1$
		\[
		\left\|\av{N}{n}T^nf_1\,S^{n^2}f_2\,S^{n^2+n}f_3\right\|_2=1.
		\]
		Furthermore, \(f_1\in L^\infty(\W_4(T))\) and
		\(f_2,f_3\in L^\infty(\W_2(S))\).
	\end{proposition}
	
	\begin{proof}
		Fix \(\alpha\in\R\setminus\Q\), and equip \(\T^5\) with Haar
		probability \(m_{\T^5}\). With all coordinates taken modulo one, define
		\[
		Tx=
		\begin{pmatrix}
			x_1+6\alpha\\
			x_2+2x_1+6\alpha\\
			x_3+3x_2+3x_1+6\alpha\\
			x_4+2x_3+3x_2+2x_1+3\alpha\\
			x_5+6x_4
		\end{pmatrix},\qquad
		Sx=
		\begin{pmatrix}
			x_1\\
			x_2+2\alpha\\
			x_3+x_1\\
			x_4+x_2+\alpha\\
			x_5+2x_3-3x_2+2x_1-3\alpha
		\end{pmatrix}.
		\]
		Their linear parts are triangular integer matrices with diagonal
		entries one, so both maps are invertible and preserve \(m_{\T^5}\).
		They commute, since substitution gives
		\[
		TSx=STx=
		\begin{pmatrix}
			x_1+6\alpha\\
			x_2+2x_1+8\alpha\\
			x_3+3x_2+4x_1+12\alpha\\
			x_4+2x_3+4x_2+4x_1+10\alpha\\
			x_5+6x_4+2x_3+3x_2+2x_1+3\alpha
		\end{pmatrix}.
		\]
		Define the circle-valued functions \(\ell_i:\T^5\to\T\) by
		\[
		\ell_1=-x_4,\qquad
		\ell_2=-x_5,\qquad
		\ell_3=x_5+3x_4-x_3,
		\qquad f_i=e(\ell_i).
		\]
		A computation of their orbit increments gives
		\[
		\begin{aligned}
			T^n\ell_1-\ell_1
			&=-2nx_3-3n^2x_2-2n^3x_1-3\alpha n^4,\\
			S^r\ell_2-\ell_2
			&=-2rx_3+3rx_2-r(r+1)x_1+3\alpha r^2,\\
			S^r\ell_3-\ell_3
			&=2rx_3+r^2x_1.
		\end{aligned}
		\]
		
		These formulas follow by induction from the definitions
		of $T$ and $S$. Substituting $r=n^2$ and $r=n^2+n$
		into the second and third formulas, respectively, gives
		
		\[
		\begin{aligned}
			(T^n\ell_1-\ell_1)
			+(S^{n^2}\ell_2-\ell_2)
			+(S^{n^2+n}\ell_3-\ell_3)
			&= \bigl(-2n-2n^2+2(n^2+n)\bigr)x_3\\
			&\quad+\bigl(-3n^2+3n^2\bigr)x_2\\
			&\quad+\bigl(-2n^3-n^2(n^2+1)
			+(n^2+n)^2\bigr)x_1\\
			&\quad+\bigl(-3n^4+3n^4\bigr)\alpha
			=0.
		\end{aligned}
		\]
		Here the coefficient of $x_1$ vanishes because
		\[
		(n^2+n)^2-n^2(n^2+1)=2n^3.
		\]
		Thus 
		\[
		\begin{aligned}
			T^nf_1\,S^{n^2}f_2\,S^{n^2+n}f_3
			=e(\ell_1+\ell_2+\ell_3)
			=e(2x_4-x_3)
		\end{aligned}
		\]
		for every $n\geq 1$. Since the right-hand side has
		modulus one, it follows that
		\[
		\left\|
		\frac1N\sum_{n=1}^N
		T^nf_1\,S^{n^2}f_2\,S^{n^2+n}f_3
		\right\|_{L^2(m_{\T^5})}=1
		\]
		for every $N\geq1$.

		For positive integers $h_1,h_2,h_3,h_4,h,k$, finite differences give
		\[
		\begin{aligned}
			\Delta_{T^{h_4}}\cdots\Delta_{T^{h_1}}f_1
			&=e(-72\alpha h_1h_2h_3h_4),\\
			\Delta_{S^k}\Delta_{S^h}f_2
			&=e\bigl(-2hk(x_1-3\alpha)\bigr),\\
			\Delta_{S^k}\Delta_{S^h}f_3
			&=e(2hkx_1).
		\end{aligned}
		\]
		The last two functions have integral zero, while
		irrationality of $\alpha$ gives
		\[
		\lim_{H\to\infty}\av{H}{h_4}
		e(-72\alpha h_1h_2h_3h_4)=0
		\]
		for fixed $h_1,h_2,h_3$. 
		By \eqref{eq:iterated-correlation},
		\[
		\|f_1\|_{U^4(T)}
		=\|f_2\|_{U^2(S)}
		=\|f_3\|_{U^2(S)}=0.
		\]
		Also,
		\[
		\Delta_T^5f_1=\Delta_S^3f_2=\Delta_S^3f_3=1.
		\]
		Hence, by the definition of phase factor,
		\[
		f_1\in L^\infty(\W_4(T)),\qquad
		f_2,f_3\in L^\infty(\W_2(S)).
		\]
		
	\end{proof}
	
	Take $T_1=T$ and $T_2=T_3=S$. This proposition shows that
$\E(f_1\mid\HK_3(T_1))=0$ and
$\E(f_i\mid\HK_1(T_i))=0$ for $i=2,3$, but the average
has $L^2$-norm one. Hence none of the three Host--Kra orders
can be lowered, even if the other two factors are unrestricted.
Since $\W_d(T_i)\subseteq\HK_d(T_i)$, the phase-factor orders
in Theorem~\ref{thm:main} are also optimal.
The summand does not depend on $n$, so the same example gives
optimality over the primes; see Section~\ref{sec:primes}.
	\section{Averages over the primes}\label{sec:primes}
	
In this section,	we prove the prime assertions of Theorems~\ref{thm:phase-reduction}
	and~\ref{thm:main}. 
	
	We need the following the comparison theorem
	of Frantzikinakis, Host, and Kra
	\cite[Proposition~3.6]{FrantzikinakisHostKra2013}.
	\begin{theorem}
		\label{thm:fhk-comparison}
		Let $T_1,\ldots,T_\ell$ be commuting invertible
		measure-preserving transformations, let
		$g_1,\ldots,g_m\in L^\infty(\mu)$, and let
		$q_{i,j}\in\Z[t]$ for $1\le i\le\ell$ and
		$1\le j\le m$.
		For $w>2$, put
		\[
		W=\prod_{\substack{p<w\\p\in\mathbb P}}p,
		\qquad
		\Lambda'_{w,b}(n)
		=\frac{\varphi(W)}{W}\Lambda'(Wn+b),
		\]
		where $\Lambda'(n)=1_{\mathbb P}(n)\log n$.
		Then
		\[
		\begin{aligned}
			\lim_{w\to\infty}\limsup_{N\to\infty}
			\max_{\substack{1\le b\le W\\(b,W)=1}}
			\bigg\|
			\frac1N\sum_{n=1}^N
			(\Lambda'_{w,b}(n)-1)
			\cdot
			\prod_{j=1}^m
			\left(\prod_{i=1}^\ell
			T_i^{q_{i,j}(Wn+b)}\right)g_j
			\bigg\|_2=0.
		\end{aligned}
		\]
	\end{theorem}
	\begin{lemma}\label{lem:prime-transfer}
		Let $T_1,\ldots,T_\ell$ be invertible commuting measure-preserving
		transformations, let $p_i\in\Z[t]$, and let $f_i\in L^\infty(\mu)$.
		Put $a(n)=\prod_{i=1}^\ell T_i^{p_i(n)}f_i$.
		If, for every $m\ge1$ and $1\le b\le m$,
		\[
		\left\|\frac1N\sum_{n=1}^N a(mn+b)\right\|_2\longrightarrow0,
		\]
		then $\|A_{N,\mathbb P}^{\mathbf p}(f_1,\ldots,f_\ell)\|_2\to0$.
	\end{lemma}
	
	\begin{proof}
		By scaling, assume $|f_i|\le1$.
		For $w>2$, put $W=\prod_{p<w,\ p\in\mathbb P}p$ and
		$\Lambda'(n)=1_{\mathbb P}(n)\log n$.
		For $1\le b\le W$ with $(b,W)=1$, define
		$\Lambda'_{w,b}(n)=\varphi(W)\Lambda'(Wn+b)/W$, where
		$\varphi(W)=\#\{1\le b\le W:(b,W)=1\}$ is Euler's totient function.
		Theorem~\ref{thm:fhk-comparison}, applied with
		polynomial $p_j$ in direction $T_j$ and zero polynomials in the other
		directions, gives
		\[
		\lim_{w\to\infty}\limsup_{N\to\infty}
		\max_{\substack{1\le b\le W\\(b,W)=1}}
		\left\|\frac1N\sum_{n=1}^N
		(\Lambda'_{w,b}(n)-1)a(Wn+b)\right\|_2=0.
		\]
		For every fixed $w$ and every $b$ with $(b,W)=1$,
		the hypothesis gives
		\[
		\frac1N\sum_{n=1}^N a(Wn+b)\longrightarrow0
		\quad\text{in }L^2(\mu).
		\]
		Since there are only finitely many such $b$,
		the convergence is uniform in $b$.
		Consequently,
		\[
		\lim_{w\to\infty}\limsup_{N\to\infty}
		\max_{\substack{1\le b\le W\\(b,W)=1}}
		\left\|\frac1N\sum_{n=1}^N
		\Lambda'_{w,b}(n)a(Wn+b)\right\|_2=0.
		\]
		Write $L_M=M^{-1}\sum_{n=1}^M\Lambda'(n)a(n)$.
		Fix $w$ and, for $M\ge W$, put $N=\lfloor M/W\rfloor$.Then 
		\[
		\sum_{b=1}^W\sum_{n=1}^N
		\Lambda'(Wn+b)a(Wn+b)
		=
		\sum_{k=W+1}^{W(N+1)}\Lambda'(k)a(k).
		\]
		Since every prime greater than $W$ is coprime to $W$ and $\Lambda'_{w,b}(n)=\varphi(W)\Lambda'(Wn+b)/W$, we have
		\[
		\begin{aligned}
			M L_M
			&=\frac{W}{\varphi(W)}
			\sum_{\substack{1\le b\le W\\(b,W)=1}}
			\sum_{n=1}^N
			\Lambda'_{w,b}(n)a(Wn+b)\\
			&\quad+\sum_{k=1}^W\Lambda'(k)a(k)
			-\sum_{k=M+1}^{W(N+1)}\Lambda'(k)a(k).
		\end{aligned}
		\]
		Since $WN\le M$ and
		$\|a(k)\|_2\le1$, we obtain
		\[
		\|L_M\|_2
		\le
		\max_{\substack{1\le b\le W\\(b,W)=1}}
		\left\|
		\frac1N\sum_{n=1}^N
		\Lambda'_{w,b}(n)a(Wn+b)
		\right\|_2
		+\frac{2W\log(M+W)}{M}.
		\]
		For fixed $w$, the last term tends to zero as
		$M\to\infty$. Taking the upper limit in $M$ and
		then letting $w\to\infty$ gives $\|L_M\|_2\to0$.

		Finally, by the  prime number theorem with $\theta(N)=\sum_{p\le N}\log p$, 
		$\|a(n)\|_2\le1$ gives
		\[
		\begin{aligned}
			\left\|\frac1{\pi(N)}\sum_{p\le N}a(p)-L_N\right\|_2
			&\le
			\left|\frac1{\pi(N)}-\frac{\log N}{N}\right|
			\sum_{p\le N}\|a(p)\|_2
			+\frac1N\sum_{p\le N}
			(\log N-\log p)\|a(p)\|_2\\
			&\le
			\left|\frac1{\pi(N)}-\frac{\log N}{N}\right|\pi(N)
			+\frac1N\sum_{p\le N}(\log N-\log p)\\
			&=
			\left|1-\frac{\pi(N)\log N}{N}\right|
			+\frac{\pi(N)\log N-\theta(N)}N\longrightarrow0.
		\end{aligned}
		\]
		This proves the lemma.
	\end{proof}

	\begin{proof}[Prime assertion of Theorem~\ref{thm:phase-reduction}]
		Let $s=s(D,\ell)$ be as in the integer assertion.
		Suppose $\E(f_j\mid\W_s(T_j))=0$ for some $j$.
		Fix $m\ge1$ and $1\le b\le m$, and put
		\[
		\widetilde p_i(t)=p_i(mt+b)-p_i(b),
		\qquad
		\widetilde f_i=T_i^{p_i(b)}f_i.
		\]
		
		The polynomials $\widetilde p_i$ have zero constant terms
		and degrees at most $D$. Moreover, for every
		nonempty $J\subseteq[\ell]$,
		$r_{\widetilde{\mathbf p}}(J)=r_{\mathbf p}(J)$, $|J|<2r_{\widetilde{\mathbf p}}(J)$.
		
		Since $\W_s(T_j)$ is $T_j$-invariant,
		$\E(\widetilde f_j\mid\W_s(T_j))=0$.
		The integer assertion therefore gives
		\[
		\frac1N\sum_{n=1}^N
		\prod_{i=1}^\ell T_i^{p_i(mn+b)}f_i
		\longrightarrow0
		\quad\text{in }L^2(\mu).
		\]
		As this holds for every $m$ and $b$,
		Lemma~\ref{lem:prime-transfer} gives the prime
		zero-limit assertion.

		Finally, let $g_i=\E(f_i\mid\W_s(T_i))$.
		Applying the prime zero-limit assertion to
		$f_i-g_i$ in each coordinate shows that
		\[
		\|A_{N,\mathbb P}^{\mathbf p}(f_1,\ldots,f_\ell)
		-A_{N,\mathbb P}^{\mathbf p}(g_1,\ldots,g_\ell)\|_2\to0.
		\]
		This proves prime characteristicity.
	\end{proof}
	
Now we show the sharp orders over the primes.
	
	\begin{proof}[Prime assertions and optimality in
		Theorem~\ref{thm:main} and Corollary~\ref{cor:orders}]
		Put $(d_1,d_2,d_3)=(4,2,2)$ and
		$a(n)=T_1^nf_1\,T_2^{n^2}f_2\,T_3^{n^2+n}f_3$.
		Suppose $\E(f_j\mid\W_{d_j}(T_j))=0$ for some $j$.
		To apply Lemma~\ref{lem:prime-transfer}, fix $m\ge1$ and
		$1\le b\le m$. Give $E=\Z/m\Z$ uniform probability $m_E$,
		let $Rz=z+1$, and define
		$\widetilde T_1=T_1\times R$,
		$\widetilde T_2=T_2\times\id_E$ and 
		$\widetilde T_3=T_3\times\id_E$.
		These transformations commute. For $z\in\Z/m\Z$,
		let
		$F_1(x,z)=f_1(x)1_{\{b+m\Z\}}(z)$, and 
		$F_i(x,z)=f_i(x)$ for $i=2,3$.

		We first check that
		$\E(F_j\mid\W_{d_j}(\widetilde T_j))=0$.
		For $\Phi\in\mathcal P_{d_j}(\widetilde T_j)$, choose $q\ge1$
		with $\Delta_{\widetilde T_j^q}^{d_j+1}\Phi=1$.
		By Lemma~\ref{lem:polynomial-phases}(i), this identity also
		holds with $qm$ in place of $q$.
		Since $\widetilde T_j^{qm}=T_j^{qm}\times\id_E$, each
		section $\Phi(\cdot,z)$ belongs to $\mathcal P_{d_j}(T_j)$.
		Since this holds for every phase generator $\Phi$,
		Lemma~\ref{lem:phasefactors}(iii) gives
		$\E(F_j\mid\W_{d_j}(\widetilde T_j))=0$.

		The integer assertion of Theorem~\ref{thm:main} now gives
		\[
		B_N(x,z):=\frac1N\sum_{n=1}^N
		1_{\{z+n=b\}}a(n)(x)\longrightarrow0
		\quad\text{in }L^2(\mu\otimes m_E).
		\]
		
		Since $E=\Z/m\Z$ has uniform measure,
		\[
		\|B_N(\cdot,0)\|_{L^2(\mu)}
		\le\sqrt m\,\|B_N\|_{L^2(\mu\otimes m_E)}
		\longrightarrow0.
		\]
		Moreover,
		\[
		\frac1N\sum_{n=1}^Na(mn+b)
		=mB_{mN}(\cdot,0)
		+\frac{a(mN+b)-a(b)}N.
		\]
		Since $\|a(n)\|_2\le\prod_{i=1}^3\|f_i\|_\infty$,
		\[
		\begin{aligned}
			\left\|\frac1N\sum_{n=1}^Na(mn+b)\right\|_2
			\le m^{3/2}\|B_{mN}\|_{L^2(\mu\otimes m_E)}
			+\frac2N\prod_{i=1}^3\|f_i\|_\infty
			\longrightarrow0.
		\end{aligned}
		\]
		As $m$ and $b$ were arbitrary,
		Lemma~\ref{lem:prime-transfer} gives the prime
		zero-limit assertion.
		Put $g_i=\E(f_i\mid\W_{d_i}(T_i))$.
		Since $\E(f_i-g_i\mid\W_{d_i}(T_i))=0$,
		we have
		\[
		\|A_{N,\mathbb P}(f_1,f_2,f_3)
		-A_{N,\mathbb P}(g_1,g_2,g_3)\|_2\to0.
		\]
		Thus the phase factors are characteristic.
		Since $\W_{d_i}(T_i)\subseteq\HK_{d_i}(T_i)$,
		the Host--Kra factors are also characteristic.
		Moreover, if $\|f_j\|_{U^{d_j+1}(T_j)}=0$
		for some $j$, the prime average tends to zero.

		Finally, take the transformations and functions
		from Proposition~\ref{prop:sharpness}, with
		$T_1=T$ and $T_2=T_3=S$. They satisfy for any $n\ge 1$,
		$T^nf_1\,S^{n^2}f_2\,S^{n^2+n}f_3
		=e(2x_4-x_3)$,
		$\|f_1\|_{U^4(T)}=0$ and
		$\|f_2\|_{U^2(S)}=\|f_3\|_{U^2(S)}=0$.
		Hence 
		$A_{N,\mathbb P}(f_1,f_2,f_3)
		=e(2x_4-x_3)$ and then
		$\|A_{N,\mathbb P}(f_1,f_2,f_3)\|_2=1$
		for every $N\ge2$. Thus all three Host--Kra orders are optimal.
		Since $\W_d\subseteq\HK_d$, the phase-factor
		orders are also optimal.
	\end{proof}
	\begin{corollary}\label{cor:phase-recurrence}
		Let $T_1,\ldots,T_\ell$ be commuting invertible
		measure-preserving transformations of a standard
		probability space $(X,\mathcal X,\mu)$, where $\ell\ge3$.
		Let $p_1,\ldots,p_\ell\in\Z[t]$ have zero constant terms, degrees at most $D$,
		satisfy \eqref{eq:rank-condition}, and
		$d_i=s(D,\ell)$ for every $i$.
		Let $A\in\mathcal X$ with $\mu(A)=\delta>0$ and 
		$\E(1_A\mid\W_{d_i}(T_i))=\delta$ for all $i\in[\ell]$.
		Then
		\[
		\frac1N\sum_{n=1}^N\mu\left(A\cap\bigcap_{i=1}^\ell
		T_i^{-p_i(n)}A\right)\longrightarrow\delta^{\ell+1},
		\qquad
		\frac1{\pi(N)}
		\sum_{\substack{p\le N\\p\in\mathbb P}}\mu\left(A\cap\bigcap_{i=1}^\ell
		T_i^{-p_i(p)}A\right)
		\longrightarrow\delta^{\ell+1}.
		\]
		Moreover, for every $0<\varepsilon<\delta^{\ell+1}$,
		the set
		$E_\varepsilon
		=\{n\ge1:\mu\left(A\cap\bigcap_{i=1}^\ell
		T_i^{-p_i(n)}A\right)>\delta^{\ell+1}-\varepsilon\}$
		satisfies
		\[
		\begin{aligned}
			\liminf_{N\to\infty}
			\frac{|E_\varepsilon\cap[1,N]|}{N}
			&\ge\frac{\varepsilon}
			{\delta-\delta^{\ell+1}+\varepsilon},\quad
			\liminf_{N\to\infty}
			\frac{|E_\varepsilon\cap\mathbb P\cap[1,N]|}{\pi(N)}
			&\ge\frac{\varepsilon}
			{\delta-\delta^{\ell+1}+\varepsilon}.
		\end{aligned}
		\]
	\end{corollary}

	\begin{proof}
	For $f_i=1_A$ ($i\in[\ell]$),	by Theorem~\ref{thm:phase-reduction},  both averages tend to $\delta^\ell$ in $L^2$.
		Multiply by $1_A$ and integrate to obtain the limit
		$\delta^\ell\mu(A)=\delta^{\ell+1}$.
		
		Let $t=\delta^{\ell+1}-\varepsilon$ and 
		$F_N=[1,N]$ or $F_N=\mathbb P\cap[1,N]$.
		Since $0\le \mu\left(A\cap\bigcap_{i=1}^\ell
		T_i^{-p_i(n)}A\right)\le\delta$ and $\mu\left(A\cap\bigcap_{i=1}^\ell
		T_i^{-p_i(n)}A\right)\le t$
		for all $n\notin E_\varepsilon$,
		\[
		\frac1{|F_N|}\sum_{n\in F_N}\mu\left(A\cap\bigcap_{i=1}^\ell
		T_i^{-p_i(n)}A\right)
		\le t+(\delta-t)
		\frac{|E_\varepsilon\cap F_N|}{|F_N|}.
		\]
		The left-hand side tends to $\delta^{\ell+1}$.
		Hence
		\[
		\liminf_{N\to\infty}
		\frac{|E_\varepsilon\cap F_N|}{|F_N|}
		\ge\frac{\delta^{\ell+1}-t}{\delta-t}
		=\frac{\varepsilon}
		{\delta-\delta^{\ell+1}+\varepsilon}.
		\]
	\end{proof}
	\begin{remark}
		For $(p_1,p_2,p_3)=(t,t^2,t^2+t)$,
		by Theorem~\ref{thm:main} one can get the same
		conclusions under the weaker assumptions
		\[
		\E(1_A\mid\W_{d_i}(T_i))=\delta,
		\qquad (d_1,d_2,d_3)=(4,2,2).
		\]
	\end{remark}
	
	\section{Consequences and a boundary example}
	\label{sec:applications}
	
	In this section, we first give bounds for the integer and prime averages
	and a necessary condition for their nonvanishing.
	We then give a Heisenberg example where phase factors imply
	convergence to zero even when the Host--Kra seminorm is nonzero.
	Finally, we take a example to show that the strict
	inequality in \eqref{eq:rank-condition} cannot
	in general be replaced by a non-strict inequality.

	\subsection{Bounds of polynomial averages}

	\begin{proposition}\label{prop:phase-control}
		Let $T_1,\ldots,T_\ell$ be commuting invertible
		measure-preserving transformations, and let
		$p_1,\ldots,p_\ell\in\Z[t]$ have zero constant terms
		and degrees at most $D$. Suppose that
		\eqref{eq:rank-condition} holds, and put $s=s(D,\ell)$
		as in Theorem~\ref{thm:phase-reduction}.
		Then for every $f_1,\ldots,f_\ell\in L^\infty(\mu)$
		with $|f_i|\le1$ for all $i\in[\ell]$, we have
		\[
		\begin{aligned}
			\limsup_{N\to\infty}
			\|A_N^{\mathbf p}(f_1,\ldots,f_\ell)\|_2
			&\le\min_{j\in[\ell]}
			\|\E(f_j\mid\W_s(T_j))\|_2,\\
			\limsup_{N\to\infty}
			\|A_{N,\mathbb P}^{\mathbf p}(f_1,\ldots,f_\ell)\|_2
			&\le\min_{j\in[\ell]}
			\|\E(f_j\mid\W_s(T_j))\|_2.
		\end{aligned}
		\]
		If either average does not tend to zero in $L^2$,
		then for every $j\in[\ell]$ there exists
		$\phi_j\in\mathcal P_s(T_j)$ such that
		\[
		\int_X f_j\overline{\phi_j}\,d\mu\ne0.
		\]

	\end{proposition}
	
	\begin{proof}
		
		Put $g_i=\E(f_i\mid\W_s(T_i))$.
		By Theorem~\ref{thm:phase-reduction}, replacing
		$f_i$ by $g_i$ does not change either average
		in the $L^2$ limit. Since $|g_i|\le1$, for every
		$j\in[\ell]$,
		\[
		\begin{aligned}
			\|A_N^{\mathbf p}(g_1,\ldots,g_\ell)\|_2
			\le\frac1N\sum_{n=1}^N
			\|T_j^{p_j(n)}g_j\|_2
			=\|g_j\|_2.
		\end{aligned}
		\]
		The same estimate holds for the prime average.
		Taking upper limits gives both inequalities.
		
		If either average does not tend to zero,
		then $g_j\ne0$ for every $j$.
		By Lemma~\ref{lem:phasefactors}(iii), there exists some
		$\phi_j\in\mathcal P_s(T_j)$ satisfies
		$\int_X f_j\overline{\phi_j}\,d\mu\ne0$.
	\end{proof}
	
	\begin{remark}
		For $(p_1,p_2,p_3)=(t,t^2,t^2+t)$,
		Theorem~\ref{thm:main} gives the same conclusions
		with $\W_s(T_j)$ replaced by
		$\W_{d_j}(T_j)$, where $(d_1,d_2,d_3)=(4,2,2)$.
		In the phase-correlation conclusion, the
		corresponding identity is
		$\Delta_{T_j^{q_j}}^{d_j+1}\phi_j=1$.
	\end{remark}

	\subsection{An example with smaller characteristic factors}
	\label{subsec:heisenberg-application}
	
	Let $X=G/\Gamma$ be the standard Heisenberg nilmanifold, where
	$G=\R^3$ has multiplication
	$(x,y,z)(x',y',z')=(x+x',y+y',z+z'+xy')$ and
	$\Gamma=\Z^3$. Let $T_2$ be an ergodic left translation and let
	$\mu$ be Haar measure. Its Kronecker factor $\mathcal K$ is
	the horizontal quotient $G/([G,G]\Gamma)$, and
	$\HK_2(T_2)=\mathcal B_X$ modulo completion.
	By \cite[Propositions~2.1(i), 2.4 and Theorem~2.2]{Frantzikinakis2008},
	every $T_2^q$ is ergodic and has order-two
	quasi-eigenfunction factor $\mathcal K$. Hence
	\[
	\W_2(T_2)=\mathcal K
	\subsetneq\HK_2(T_2)=\mathcal B_X.
	\]
	
	Choose a nonzero $g\in L^\infty(\mu)$ with $\E(g\mid\mathcal K)=0$
	Then $\|g\|_{U^3(T_2)}>0$, since the $U^3(T_2)$ null space is
	orthogonal to $L^2(\HK_2(T_2))=L^2(\mu)$.
	However, for any  $f_1,f_3\in L^\infty(\mu)$ and any $T_1,T_3$
	commuting with each other and with $T_2$, by Theorem~\ref{thm:main}, we have
	\[
	\|A_N(f_1,g,f_3)\|_2\longrightarrow0,
	\qquad
	\|A_{N,\mathbb P}(f_1,g,f_3)\|_2\longrightarrow0.
	\]
	Hence these averages can tend to zero even if the Host--Kra
	seminorm is nonzero.

	\subsection{A boundary example for the rank condition}
	\label{subsec:rank-boundary}
	
    In general,	the strict inequality in \eqref{eq:rank-condition}
	cannot be replaced by $|J|\le2r_{\mathbf p}(J)$. The following example shows that,
	even when all transformations are equal, no
	finite-order phase factor need be characteristic.

	\begin{proposition}\label{prop:rank-boundary}
		Let $p_j(t)=t+jt^2$ for $j=0,1,2,3$. Then
		there exist a totally ergodic nilrotation $T$ on a Heisenberg nilmanifold $X$ and functions
		$f_0,f_1,f_2,f_3\in L^\infty(\mu)$ with $|f_j|\le1$
		such that, for every $ s\ge1,0\le j\le3$
		\[
		\E(f_j\mid\W_s(T))=0,
		\]
		whereas
		\[
		\left\|
		\frac1N\sum_{n=1}^N
		\prod_{j=0}^3T^{p_j(n)}f_j-\frac1{12}
		\right\|_2\longrightarrow0.
		\]
		Moreover, for every nonempty
		$J\subseteq\{0,1,2,3\}$,
		\[
		r_{\mathbf p}(J)=\min\{|J|,2\},
		\qquad |J|\le2r_{\mathbf p}(J),
		\]
		with equality in the second inequality only when
		$J=\{0,1,2,3\}$.
	\end{proposition}

	\begin{proof}
		
		Take the nilrotation $T:u\mapsto au$ on the Heisenberg nilmanifold from
		Subsection~\ref{subsec:heisenberg-application},
		where $a=(\alpha,\beta,\gamma)$.
		Then $T$ is totally ergodic,
		$\beta$ is irrational, and $\W_2(T)$ is equal to its Kronecker factor $K$.
		Since the eigenfunctions of $T$ generate $\mathcal K$,
		$\mathcal K\subseteq\W_1(T)
		\subseteq\W_2(T)=\mathcal K$.
		Thus $\W_1(T)=\W_2(T)=\mathcal K$, and
		Lemma~\ref{lem:phasefactors}(ii) gives
		$\W_1(T^q)=\W_2(T^q)=\mathcal K$ for all
		$q\ge1$.

		We first show this equality holds for all finite order.
		Fix $R=T^q$ and  $|\phi|=1$ with
		$\Delta_R^{s+1}\phi=1$.
		Let $r\ge0$ be least with $\Delta_R^{r+1}\phi=1$.
		If $r\ge2$, then $\psi=\Delta_R^{r-2}\phi$ belongs to
		$\mathcal P_2(R)\subseteq L^\infty(\mathcal K)$, and
		by ergodicity, $\lambda=\Delta_R^2\psi$ is constant, and the minimality of $r$ gives $\lambda\ne1$.
		Lemma~\ref{lem:phaseprojection}, with $d=1$ and $L=3$, gives
		$\psi=\E(\psi\mid\W_1(R))=1_{\{\lambda\in\roots\}}\psi$.
		Thus $\lambda$ is a root of unity. But
		$R(\Delta_R\psi)=\lambda\Delta_R\psi$, contradicting total
		ergodicity. Hence $r\le1$, so $\phi\in\mathcal P_1(R)$.
		This proves $\W_s(T)=\mathcal K$ for all $s\ge1$.
		
		Use the fundamental domain $[0,1)^3$, with coordinates $(x,y,z)$
		and Lebesgue measure. Put $I=[0,1/2)$ and
		\[
		(c_0,c_1,c_2,c_3)=(-1,3,-3,1),\qquad
		f_j(x,y,z)=1_I(y)e(c_jz).
		\]
		Since $\mathcal K$ is the horizontal factor and
		$c_j\ne0$ for every $j$, we have
		\[
		\E(f_j\mid\mathcal K)(x,y)
		=1_I(y)\int_0^1 e(c_jz)\,dz=0.
		\]
		Denote
		$C_N=\frac1N\sum_{n=1}^N
		\prod_{j=0}^3 T^{n+jn^2}f_j$.
		We shall prove that
		$\left\|C_N-\frac1{12}\right\|_2\longrightarrow0$.

		Fix $n\ge1$ and put $t_j=n+jn^2$ for $0\le j\le3$.
		Recall that for any $t\in\Z$
		\[
		a^t(x,y,z)=
		\left(x+t\alpha,\ y+t\beta,\
		z+t\gamma+\binom t2\alpha\beta+t\alpha y\right).
		\]
		Put $Y_j=y+\beta t_j$ and $k_j=\lfloor Y_j\rfloor$.
		After reducing $a^{t_j}(x,y,z)$ modulo $\Gamma$,
		its second coordinate is $\{Y_j\}$, and its third
		coordinate is
		$q_j:=z+\gamma t_j+\binom{t_j}{2}\alpha\beta
		+\alpha t_jy-(x+\alpha t_j)k_j$ modulo one.

		If some $\{Y_j\}\notin I$, then $\prod_{j=0}^3T^{t_j}f_j(x,y,z)
		=
		\prod_{j=0}^3 1_I(\{Y_j\})
		\cdot e\left(\sum_{j=0}^3c_jq_j\right)=0.$
		Now we consider the case of all $\{Y_j\}$ belong to $[0,1/2)$.
		Since $Y_j$ has degree one in $j$, $
		Y_{j+2}-2Y_{j+1}+Y_j=0$. Then we have
		$k_{j+2}-2k_{j+1}+k_j
		=-(\{Y_{j+2}\}-2\{Y_{j+1}\}+\{Y_j\})=0$
		for $j=0,1$.
		Thus $k_j=k_0+j(k_1-k_0)$.
		
		Since $t_j$ and $k_j$ have degree at most one
		in $j$, $q_j$ has degree at most
		two in $j$. By
		$(c_0,c_1,c_2,c_3)=(-1,3,-3,1)$, $\sum_{j=0}^3c_jq_j
		=-q_0+3q_1-3q_2+q_3=0$.
		Consequently,
		\[
		\begin{aligned}
			C_N(x,y,z)
			=\frac1N\sum_{n=1}^N
			\prod_{j=0}^3 1_I(\{y+n\beta+jn^2\beta\})
			=\frac1N\sum_{n=1}^N F(y+n\beta,n^2\beta),
		\end{aligned}
		\]
		where $F(s,t)=\prod_{j=0}^3 1_I(\{s+jt\})$.

		Take $0\le s<1$ and $-1/2\le t<1/2$.
		If $F(s,t)=1$, then $\{s+jt\}\in[0,1/2)$
		for $0\le j\le3$. For $j=0,1,2$,
		\[
		\begin{aligned}
			\lfloor s+(j+1)t\rfloor-\lfloor s+jt\rfloor
			=t-\{s+(j+1)t\}+\{s+jt\}
			\in(-1,1).
		\end{aligned}
		\]
		This difference is an integer, so it is zero.
		Since $\lfloor s\rfloor=0$, we obtain
		$\lfloor s+jt\rfloor=0$ for all $j$.
		Hence $F(s,t)=1$  when
		$0\le s+jt<1/2$ for $0\le j\le3$. Then
		\[
		\begin{aligned}
			F(s,t)=1
			&\iff 0\le s<\frac12,\quad
			-\frac{s}{3}\le t<\frac{1/2-s}{3},\\
			\int_{\T^2}F(s,t)\,ds\,dt
			&=\int_0^{1/2}\frac16\,ds=\frac1{12}.
		\end{aligned}
		\]
		
		Since $\beta$ is irrational, Weyl's theorem
		\cite{Weyl1916} gives equidistribution of
		$(y+n\beta,n^2\beta)$ in $\T^2$ for every $y$.
		The function $F$ is Riemann integrable,
		since its discontinuities lie on the lines
		$s+jt\equiv0,1/2\pmod1$.
		Thus $C_N(x,y,z)\to1/12$.
		Since $0\le C_N\le1$, bounded convergence gives
		\[
		\left\|C_N-\frac1{12}\right\|_2\longrightarrow0.
		\]

		Finally, since $p_j(t)=t+jt^2$, any two distinct
		polynomials $p_i$ and $p_j$ are linearly independent
		over $\Q$.
		Hence $r_{\mathbf p}(J)=\min\{|J|,2\}$,
		as required.
		
	\end{proof}

\end{document}